\documentclass[smallextended,natbib]{svjour3}
\smartqed  
\spnewtheorem{algorithm}{Algorithm}{\bf}{\rm}
\usepackage{amsmath,amssymb,enumerate}
\usepackage{graphicx}
\usepackage{mathptmx,bm}\usepackage[mathcal]{euscript}
\usepackage{helvet,courier}

\usepackage{subfigure}\renewcommand{\subfigtopskip}{0pt}\renewcommand{\subfigcapskip}{0pt}
\usepackage{mathtools}\mathtoolsset{showonlyrefs=false}
\usepackage{booktabs}
\newcommand{\refalg}[1]{Algorithm~\ref{#1}}
\newcommand{\refprop}[1]{Proposition~\ref{#1}}

\newcommand{\reffig}[1]{Figure~\ref{#1}}
\newcommand{\reftab}[1]{Table~\ref{#1}}

\newcommand{\finbox}{\nolinebreak\hfill{\small $\blacksquare$}}

\newcounter{alnum}
\newenvironment{algstep}
{\begin{list}{{\upshape\bf Step~\arabic{alnum}:}}%
{\usecounter{alnum}%
\setlength{\leftmargin}{18mm}%
\setlength{\labelwidth}{20mm}%
}}{\end{list}}

\newenvironment{Problem}{\begin{array}.{*{20}{l}}\}}{\end{array}}

\newcommand{\Min}{\mathop{\mathrm{minimize}}}

\newcommand{\ST}{\mathop{\mathrm{subject~to}}}
\newcommand{\sign}{\mathop{\mathrm{sgn}}\nolimits}

\newcommand{\diag}{\mathop{\mathrm{diag}}\nolimits}

\newcommand{\argmin}{\operatornamewithlimits{\mathrm{arg\,min}}}
\newcommand{\argmax}{\operatornamewithlimits{\mathrm{arg\,max}}}

\renewcommand{\Re}{\ensuremath{\mathbb{R}}}

\newcommand{\bi}[1]{\ensuremath{\boldsymbol{#1}}}
\newcommand{\ti}[1]{\ensuremath{#1}}

\newcommand{\rr}[1]{\ensuremath{\mathrm{#1}}}
\newcommand{\bs}[1]{\ensuremath{\boldsymbol{\mathsf{#1}}}}

\newcommand{\pdif}[2]{\frac{\partial #1}{\partial #2}}

\newcommand{\LC}{\ensuremath{\mathcal{L}}}

\journalname{Optimization and Engineering}

\begin{document}

\title{A fast first-order optimization approach 
  to elastoplastic analysis of skeletal structures 
  }

\titlerunning{A fast first-order optimization approach to elastoplastic analysis 
  }

\author{Yoshihiro Kanno
  }

\institute{
  Yoshihiro Kanno \at
  Address: 
  Laboratory for Future Interdisciplinary Research of Science and Technology, 
  Institute of Innovative Research, 
  Tokyo Institute of Technology, 
  Nagatsuta 4259, Yokohama 226-8503, Japan. \\
  E-mail: \texttt{kanno.y.af@m.titech.ac.jp} \\
  Phone: +81-45-924-5364 \\
  Fax: +81-45-924-5977
  }

\date{Received: date / Accepted: date}

\maketitle

\begin{abstract}
  It is classical that, when the small deformation is assumed, 
  the incremental analysis 
  problem of an elastoplastic structure with a piecewise-linear yield 
  condition and a linear strain hardening model can be formulated as a 
  convex quadratic programming problem. 
  Alternatively, 
  this paper presents a different formulation, an unconstrained 
  nonsmooth convex optimization problem, and proposes to solve it with 
  an accelerated gradient-like method. 
  Specifically, we adopt an accelerated proximal gradient method, that 
  has been developed for a regularized least squares problem. 
  Numerical experiments show that the presented algorithm is effective 
  for large-scale elastoplastic analysis. 
  Also, a simple warm-start strategy can speed up the algorithm when 
  the path-dependent incremental analysis is carried out. 
\end{abstract}

\section{Introduction}

It has been diversely recognized that the elastoplastic incremental 
analysis of solids and structures is very linked to theory and 
algorithms of optimization; see, e.g., \citet{MM82} for survey. 
If the small deformation is assumed, the incremental problem of an 
elastoplastic truss can be formulated as a linear complementarity 
problem (LCP) \citep{dDM76,Smi78,Kaneko79,Kaneko80,WTl90,TlX01,TTl07,TTl08}. 
It is well known that, if the hardening modulus is nonnegative (i.e., if 
the strain softening is not considered), then this LCP can be recast as 
a (convex) quadratic programming (QP) 
problem \citep{Mai68,Mai70,CM70,GFDdC79,SP10}. 
We can solve a QP problem efficiently with a primal-dual 
interior-point method \citep{AL12}. 

This paper attempts to shed new light on this classical problem in 
computational plasticity from perspective of a recently developed 
branch of numerical optimization. 
Namely, in this paper we examine a simpler 
gradient-based algorithm with acceleration. 
Instead of QP and LCP, we formulate the incremental problem as an 
unconstrained nonsmooth convex optimization problem. 

Recently, accelerated, or ``optimal'' \citep{Nes04}, first-order methods 
have received considerable attention, particularly for solving 
large-scale optimization problems arising in image processing, 
regression, etc.; see, e.g., 
\citet{BT09}, \citet{GOSB14}, and \citet{OdC15}. 
Such a method converges in the objective value with rate $O(1/k^{2})$, 
where $k$ is the iteration counter. 
Also, since it is basically a gradient-like method, the computation 
at each iteration is very cheap. 

In this paper we show that the incremental problem of an elastoplastic 
truss can be solved efficiently with an accelerated first-order method. 
Specifically, we adopt an accelerated proximal gradient method 
\citep{BT09,CEg14,PB14,OdC15}. 
Computational effort at each iteration of the presented algorithm is 
dominated by cheap matrix-vector multiplications. 
When we adopt a variant of the Newton--Raphson method for elastoplastic 
analysis, a major difficulty is to find a consistent tangent stiffness 
matrix, because one cannot know {\em a priori\/} each structural element 
will undergo plastic loading or elastic unloading. 
As a common attribute among optimization approaches to elastoplastic 
analysis, the presented algorithm does not use the tangent stiffness 
matrix, and automatically detects whether each member will 
undergo plastic loading or elastic unloading. 
Moreover, unlike other popular efficient optimization algorithms, the 
presented algorithm does not use a linear-equations solver. 

In the course of path-dependent quasi-static analysis, we solve the 
incremental problem repeatedly with varying the load parameter and 
updating the state variables. 
This means that we solve a series of closely related problems. 
Since the presented approach is based upon the unconstrained 
optimization formulation, it might possibly employ a simple warm-start 
strategy that uses the solution at the previous loading step as the 
initial solution for the present loading step. 
The effect of this warm-start strategy will be investigated through 
numerical experiments. 
In contrast, interior-point methods usually require some specific 
techniques for warm start; see, e.g., 
\citet{Mit01}, \citet{BS07}, \citet{JY08}, and \citet{YK12}. 



The paper is organized as follows. 
Section~\ref{sec.fundamental} summarizes fundamentals of the incremental 
analysis of an elastoplastic truss. 
As the major contribution, section~\ref{sec.proximal} presents an 
accelerated proximal gradient method for solving the incremental problem. 
Section~\ref{sec.kinematic} extends the method to a mixed model of 
isotropic hardening and kinematic hardening. 
An extension to a piecewise-linear hardening model is presented in 
section~\ref{sec.piecewise}. 
Section~\ref{sec.ex} reports the results of numerical experiments. 
We conclude in section~\ref{sec.conclude}. 


A few words regarding notation. 
We use ${}^{\top}$ to denote the transpose of a vector or a matrix. 
For simplicity, we often write the $(n+m)$-dimensional column vector 
$(\bi{x}^{\top} , \bi{y}^{\top})^{\top}$ consisting of 
$\bi{x} \in \Re^{n}$ and $\bi{y} \in \Re^{m}$ as 
$(\bi{x}, \bi{y})$. 
The Euclidean norm and the $\ell_{\infty}$-norm 
of $\bi{x}=(x_{i})\in\Re^{n}$ are denoted by 
$\| \bi{x} \| = \sqrt{\bi{x}^{\top} \bi{x}}$ and 
$\| \bi{x} \|_{\infty}=\max \{ |x_{1}|,\dots,|x_{n}| \}$, respectively. 
For a closed convex function $f:\Re^{n} \to \Re$, we define the proximal 
mapping of $f$ by 
\begin{align*}
  \bs{prox}_{f}(\bi{x})
  = \argmin_{\bi{z}} \Bigl\{
  f(\bi{z}) + \frac{1}{2} \| \bi{z} - \bi{x} \|^{2} 
  \Bigr\} . 
\end{align*}
We use $\partial f(\bi{x}) \subseteq \Re^{n}$ to denote the 
subdifferential of $f$ at a point $\bi{x} \in \Re^{n}$. 
The signum function is denoted by $\sign$, i.e., 
\begin{align*}
  \sign(s) &= 
  \begin{dcases*}
    1  & if $s > 0$, \\
    0  & if $s = 0$, \\
    -1 & if $s < 0$ 
  \end{dcases*}
\end{align*}
for $s \in \Re$. 
For vectors $\bi{x}=(x_{i}) \in \Re^{n}$ and 
$\bi{z}=(z_{i}) \in \Re^{n}$, we use 
$|\bi{x}| \in \Re^{n}$, $\sign(\bi{x}) \in \Re^{n}$, 
and $\max\{ \bi{x},\bi{z} \} \in \Re^{n}$ to denote 
\begin{align*}
  |\bi{x}| 
  &= (|x_{1}|, \dots, |x_{n}|)^{\top} , \\
  \sign(\bi{x}) 
  &= (\sign(x_{1}), \dots, \sign(x_{n}))^{\top} , \\
  \max \{\bi{x}, \bi{z}\} 
  &= (\max \{ x_{1},z_{1} \}, \dots, \max \{ x_{n},z_{n} \})^{\top} . 
\end{align*}
We use $\diag(\bi{x})$ to denote a diagonal matrix, the vector of 
diagonal components of which is $\bi{x}$.

\section{Fundamentals of elastoplastic analysis}
\label{sec.fundamental}

In this section we recall the quasi-static analysis of an elastoplastic 
truss and formulate the incremental problem; 
see, e.g., \citet{SH98}, \citet{dSnPO08}, and \citet{HR13} 
for fundamentals of computational plasticity. 

Consider an elastoplastic truss in the two- or three-dimensional space. 
Throughout the paper we assume small deformation. 
We use $m$ and $d$ to denote the number of members and the number of 
degrees of freedom of the displacements, respectively. 
In this section we consider an isotropic hardening model; see 
section~\ref{sec.kinematic} for kinematic hardening. 

Suppose that change in the external forces applied to the truss occurs 
quite slowly. 
Then the inertial term of the equations of motion becomes negligibly 
small. 
Therefore, we omit the inertial term when we construct the governing 
equations. 
The structural behavior modeled in this manner is neither static nor 
dynamic, and is referred to as quasi-static. 
We use the term ``time'' to stand for a parameter with respect to which 
the evolution process of the quasi-static behavior is described. 
This parameter, sometimes called pseudo-time, needs not correspond to 
the actual time, because the quasi-static behavior differs from the 
dynamic one. 
Suppose that we shall investigate quasi-static response of the truss 
within the time interval $[0, T]$. 
This time interval is  subdivided into finitely many intervals. 
For a specific subinterval, denoted $[t, t+\Delta t]$, 
the response at time $t+\Delta t$ is found by applying the 
standard backward (or fully implicit) Euler scheme. 

Let $\bi{u} \in \Re^{d}$ and $\bi{f} \in \Re^{d}$ denote the 
vector of displacements and the vector of external forces, respectively. 
It should be clear that these are values at time $t+\Delta t$. 
We attempt to compute $\bi{u}$ when $\bi{f}$ is specified. 
With the superscript ${}^{(t)}$ we denote the values of variables 
at time $t$ (e.g., $\bi{u}^{(t)}$ for the displacement), 
and with the prefix $\Delta$ we denote the increments between time $t$ 
and $t+\Delta t$ (e.g., $\Delta\bi{u}$ for the incremental displacement). 
The values at time $t$ are supposed to be known, and hence $\bi{u}$ is 
obtained by finding $\Delta\bi{u}$. 

Let $c_{i}$ denote the elongation of member $i$ $(i=1,\dots,m)$. 
The compatibility relation between the incremental member elongation and 
the incremental displacements can be written in the form 
\begin{align}
  \Delta c_{i} = \bi{b}_{i}^{\top} \Delta\bi{u} , 
  \label{eq.fundamental.truss.2}
\end{align}
where $\bi{b}_{i} \in \Re^{d}$ is a constant vector. 
We decompose $\Delta c_{i}$ additively as 
\begin{align}
  \Delta c_{i} 
  = \Delta c_{\rr{e}i} + \Delta c_{\rr{p}i} , 
  \label{eq.fundamental.truss.1}
\end{align}
where $\Delta c_{\rr{e}i}$ and $\Delta c_{\rr{p}i}$ are the elastic and 
plastic parts, respectively. 

Let $q_{i}$ denote the axial force of member $i$ at time $t + \Delta t$, 
which is written as 
\begin{align}
  q_{i} 
  = q^{(t)}_{i} + \Delta q_{i} . 
  \label{eq.fundamental.truss.3}
\end{align}
The constitutive law is written in terms of the increments as 
\begin{align}
  \Delta q_{i} = k_{i} \Delta c_{\rr{e}i} , 
  \label{eq.fundamental.truss.4}
\end{align}
where the elongation stiffness, $k_{i} > 0$, is assumed to be constant. 
Specifically, we have $k_{i} = E a_{i}/l_{i}$, where $E$ is the Young 
modulus, $a_{i}$ is the member cross-sectional area, and $l_{i}$ is the 
undeformed member length. 
The force-balance equation between the external forces and member 
axial forces at time $t+\Delta t$ can be written as 
\begin{align}
  \sum_{i=1}^{m} q_{i} \bi{b}_{i} = \bi{f} . 
  \label{eq.fundamental.truss.5}
\end{align}

The yield condition is given by 
\begin{align}
  |q_{i}| - R_{i} \le 0 , 
  \label{eq.fundamental.yield.1}
\end{align}
where $R_{i}$ corresponds to the magnitude of yield axial force. 
Define $\Delta\gamma_{i}$ by 
\begin{align}
  \Delta\gamma_{i}
  &= |\Delta c_{\rr{p}i}| , 
  \label{eq.def.variable.gamma}
\end{align}
which is the integration of the plastic multiplier between time $t$ and 
$t+\Delta t$. 
Under the hypothesis of linear isotropic hardening, the evolution 
of $R_{i}$ is written in the form 
\begin{align}
  R_{i} 
  &= R^{(t)}_{i} + h_{i} \Delta\gamma_{i} , 
  \label{eq.fundamental.truss.7}
\end{align}
where $h_{i} > 0$ is a constant called the isotropic hardening modulus. 
As usual, we postulate the normality flow rule, that is written as 
\begin{subequations}\label{eq.normality.low.1}
  \begin{alignat}{3}
    & q_{i} = R_{i} 
    &{\quad}& \Rightarrow {\quad} \Delta c_{\rr{p}i} \ge 0 , \\
    & q_{i} = -R_{i} 
    &{\quad}& \Rightarrow {\quad} \Delta c_{\rr{p}i} \le 0 , \\
    & |q_{i}| < R_{i} 
    &{\quad}& \Rightarrow {\quad} \Delta c_{\rr{p}i} = 0  . 
  \end{alignat}
\end{subequations}
In other words, $q_{i}$ should satisfy 
\begin{align*}
  q_{i} \in \argmax_{\hat{q}_{i}} 
  \{ \hat{q}_{i} \Delta c_{\rr{p}i} \mid 
  |\hat{q}_{i}| \le R_{i} \} , 
\end{align*}
which is called the principle of maximum plastic work. 
Here, the objective function is the plastic work due to the incremental 
plastic elongation, and the constraint is the yield condition. 
Namely, this principle states that $q_{i}$ corresponding to 
$\Delta c_{\rr{p}i}$ is the one maximizing the plastic work among the 
axial forces satisfying \eqref{eq.fundamental.yield.1}. 
This is the most fundamental and widely accepted hypothesis in the 
plasticity theory. 

\begin{figure}[tp]
  \centering
  \subfigure[]{
  \label{fig.plastic_law_1}
  \includegraphics[scale=0.65]{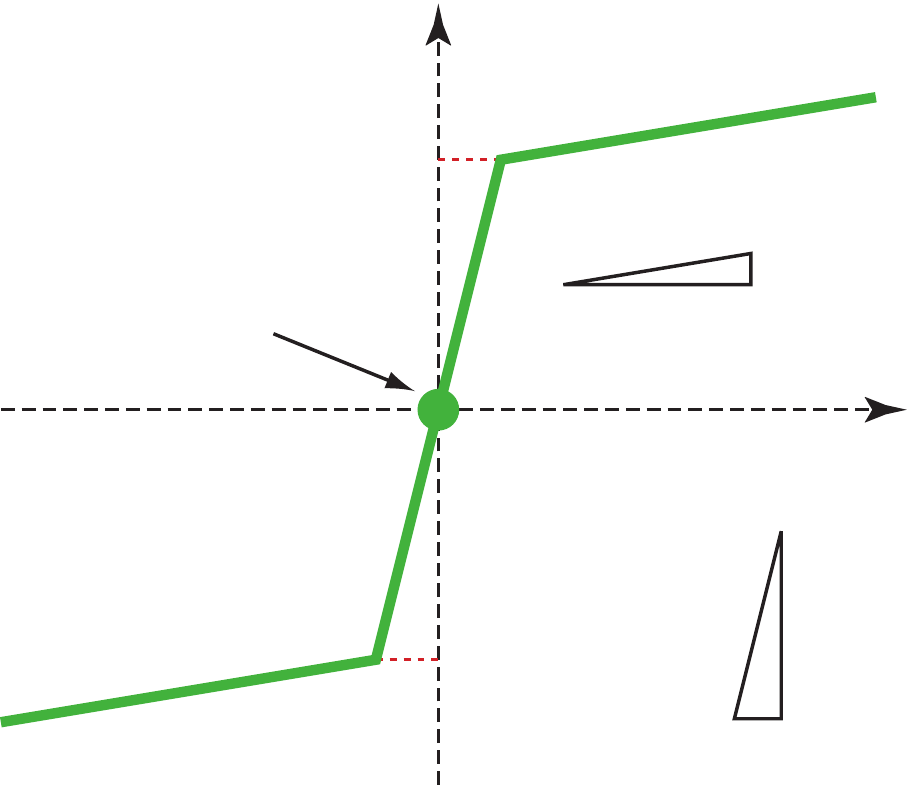}
  \begin{picture}(0,0)
    \put(-180,-10){
    \put(79,150){{\small $q_{i}$}}
    \put(140,70){{\small $c^{(0)}_{i} + \Delta c_{i}$}}
    \put(92,32){{\small $-R^{(0)}_{i}$}}
    \put(73,124){{\small $R^{(0)}_{i}$}}
    \put(148,14){{\small $1$}}
    \put(156,36){{\small $k_{i}$}}
    \put(132,96){{\small $1$}}
    \put(149,105){{\small $\displaystyle\frac{k_{i}h_{i}}{k_{i}+h_{i}}$}}
    \put(21,94){{\small $(c^{(0)}_{i},q^{(0)}_{i})$}}
    }
  \end{picture}
  }
  \hfill
  \subfigure[]{
  \label{fig.plastic_law_2}
  \includegraphics[scale=0.65]{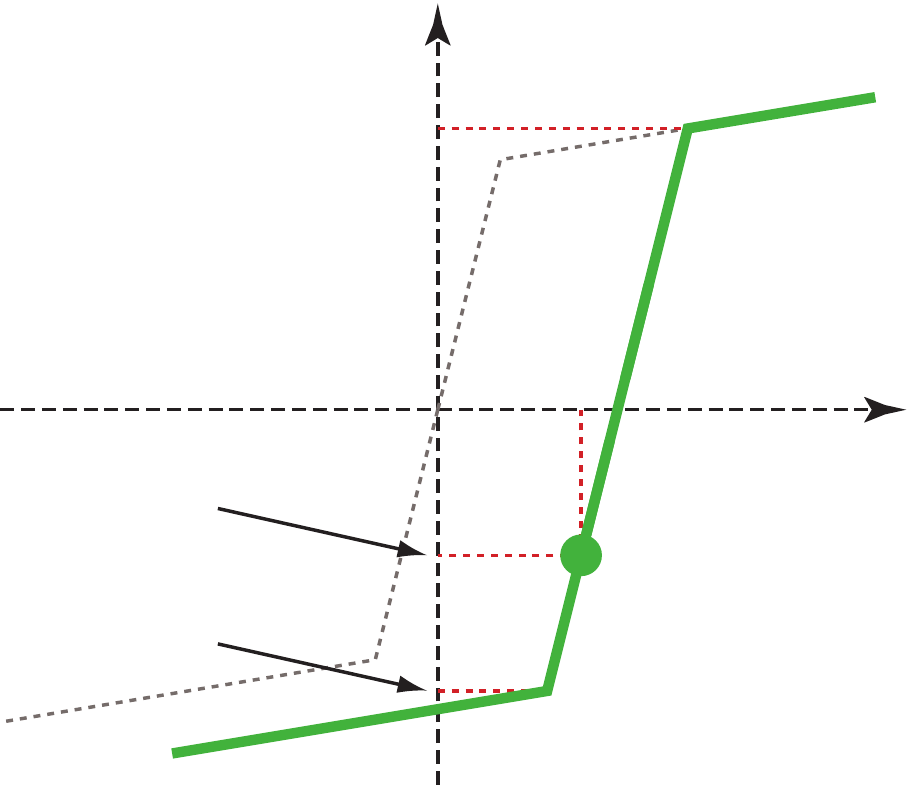}
  \begin{picture}(0,0)
    \put(-180,-10){
    \put(79,150){{\small $q_{i}$}}
    \put(140,70){{\small $c^{(t)}_{i} + \Delta c_{i}$}}
    \put(29,36){{\small $-R^{(t)}_{i}$}}
    \put(73,129){{\small $R^{(t)}_{i}$}}
    \put(110,85.5){{\small $c^{(t)}_{i}$}}
    \put(37,61){{\small $q^{(t)}_{i}$}}
    }
  \end{picture}
  }
  \caption{The constitutive law at 
  \subref{fig.plastic_law_1} time $t=0$; and 
  \subref{fig.plastic_law_2} time $t=t$, after plastic deformation has 
  taken place. }
  \label{fig.plastic_law}
\end{figure}
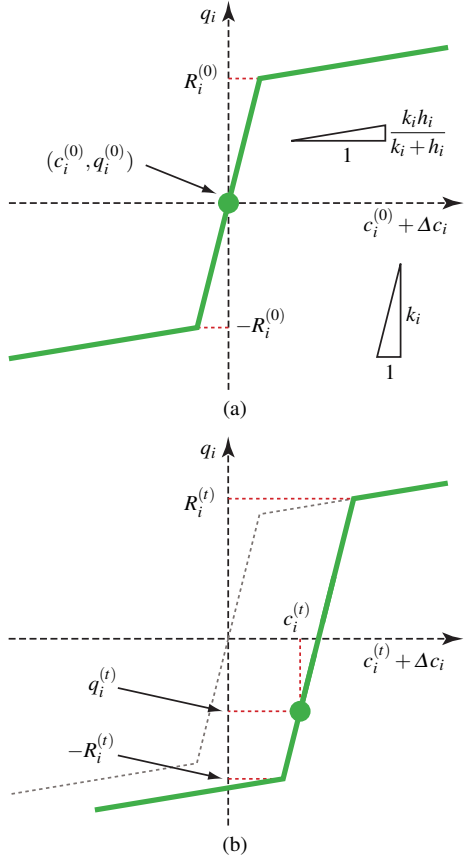

\reffig{fig.plastic_law} shows the relationship between $q_{i}$ and 
$\Delta c_{i}$ defined by \eqref{eq.fundamental.truss.1}, 
\eqref{eq.fundamental.truss.3}, \eqref{eq.fundamental.truss.4}, 
\eqref{eq.fundamental.yield.1}, \eqref{eq.def.variable.gamma}, 
\eqref{eq.fundamental.truss.7}, and \eqref{eq.normality.low.1}. 
At time $0$, we usually suppose that $c^{(0)}_{i}=q^{(0)}_{i}=0$ as 
depicted with a filled circle in \reffig{fig.plastic_law_1}. 
Member $i$ undergoes plastic deformation between time $0$ and $\Delta t$ 
if $q_{i}$ becomes greater than $R^{(0)}_{i}$. 
If this is the case, then $R^{(t)}_{i} > R^{(0)}_{i}$ as shown in 
\reffig{fig.plastic_law_2}. 

The following is a key to our formulation. 
\begin{proposition}\label{prop.complementarity}
  Assume $R_{i} > 0$. 
  Then $R_{i}$, $q_{i}$, $\Delta\gamma_{i}$, and $\Delta c_{\rr{p}i}$ 
  satisfy \eqref{eq.fundamental.yield.1}, \eqref{eq.def.variable.gamma}, 
  and \eqref{eq.normality.low.1} if and only if they satisfy 
  \begin{align}
    R_{i} \ge |q_{i}| , 
    \quad
    \Delta\gamma_{i} \ge |{-}\Delta c_{\rr{p}i}| , 
    \quad
    \begin{bmatrix}
      R_{i} \\ q_{i} \\
    \end{bmatrix}
    ^{\top}
    \begin{bmatrix}
      \Delta\gamma_{i} \\ -\Delta c_{\rr{p}i} \\
    \end{bmatrix}
    = 0 . 
    \label{eq.normality.low.2}
  \end{align}
\end{proposition}
We omit the proof; see \citet[Proposition~3]{YK12}. 
The two inequalities in \eqref{eq.normality.low.2} can be viewed as the 
second-order cone constraints in the two-dimensional space. 
The equation in \eqref{eq.normality.low.2} is then understood as a 
complementarity condition over the second-order cones; 
see, e.g., \citet{BtN01}, \citet{AL12}, and \citet{Kan11} 
for fundamentals of second-order cone constraints and complementarity 
conditions. 

We always have $R^{(0)}_{i} > 0$ and we assume $h_{i} > 0$. 
This and \eqref{eq.fundamental.truss.7} imply that the assumption made 
in \refprop{prop.complementarity} is satisfied at any time $t$. 
Accordingly, from \eqref{eq.fundamental.truss.2}, 
\eqref{eq.fundamental.truss.1}, \eqref{eq.fundamental.truss.3}, 
\eqref{eq.fundamental.truss.4}, \eqref{eq.fundamental.truss.5}, 
\eqref{eq.fundamental.truss.7}, and 
\refprop{prop.complementarity}, 
the incremental problem to be solved can be formulated as 
\begin{subequations}\label{eq.complementarity.form}
  \begin{alignat}{2}
    & \Delta c_{\rr{e}i} + \Delta c_{\rr{p}i} 
    = \bi{b}_{i}^{\top} \Delta\bi{u} , 
    &{\quad}& i=1,\dots,m, 
    \label{eq.complementarity.form.1} \\
    & q_{i} = q^{(t)}_{i} + k_{i} \Delta c_{\rr{e}i} , 
    &{\quad}& i=1,\dots,m, 
    \label{eq.complementarity.form.3} \\
    & \sum_{i=1}^{m} q_{i} \bi{b}_{i} = \bi{f} , 
    \label{eq.complementarity.form.4} \\
    & R^{(t)}_{i} + h_{i} \Delta\gamma_{i} \ge |q_{i}| , 
    \quad
    \Delta\gamma_{i} \ge |\Delta c_{\rr{p}i}| , 
    \notag\\
    &{\qquad} 
    \begin{bmatrix}
      R^{(t)}_{i} + h_{i} \Delta\gamma_{i} \\ q_{i} \\
    \end{bmatrix}
    ^{\top}
    \begin{bmatrix}
      \Delta\gamma_{i} \\ {-}\Delta c_{\rr{p}i} \\
    \end{bmatrix}
    = 0 , 
    &{\quad}& i=1,\dots,m. 
    \label{eq.complementarity.form.5}
  \end{alignat}
\end{subequations}
It should be clear in \eqref{eq.complementarity.form} that 
$\Delta\bi{u}$, $\Delta c_{\rr{e}i}$, 
$\Delta c_{\rr{p}i}$, $\Delta\gamma_{i}$, $q_{i}$ $(i=1,\dots,m)$ are 
variables to be found. 
Problem \eqref{eq.complementarity.form} is a second-order cone linear 
complementarity problem (SOCLCP). 
It is known that SOCLCP and the second-order cone programming (SOCP) 
have diverse applications in applied mechanics, including frictional 
contact \citep{Kan11,KMPdc06}, cable networks \citep{KOI02}, 
and elastoplastic continua \citep{BMP05,Mak06,KLS07,KL12,YK12}. 

\begin{remark}
  It is known that the incremental problem of an elastoplastic 
  truss can be formulated as a linear complementarity 
  problem (LCP); see, e.g., \citet{dDM76}, \citet{Kaneko79}, 
  \citet{Smi78}, \citet{TlX01}, and \citet{WTl90}. 
  Indeed, it is possible to recast \eqref{eq.complementarity.form} as an 
  LCP by splitting variables as 
  $\Delta c_{\rr{p}i}=\Delta c_{\rr{p}i}^{+}-\Delta c_{\rr{p}i}^{-}$ and 
  $\Delta \gamma_{i}=\Delta c_{\rr{p}i}^{+}+\Delta c_{\rr{p}i}^{-}$ 
  with $\Delta c_{\rr{p}i}^{+} \ge 0$ and $\Delta c_{\rr{p}i}^{-} \ge 0$ 
  and replacing the complementarity conditions by 
  $\Delta c_{\rr{p}i}^{+} (q_{i}-R_{i}) = 0$ and 
  $\Delta c_{\rr{p}i}^{-} (q_{i}+R_{i}) = 0$. 
  The resulting LCP has $2m$ complementarity conditions, while SOCLCP 
  \eqref{eq.complementarity.form} has $m$ complementarity conditions 
  (over the second-order cones). 
  It has been well recognized that this LCP can be recast as 
  (convex) quadratic programming (QP); see, 
  e.g., \citet{Mai68,Mai70}, \citet{CM70}, and \citet{GFDdC79}. 
  In contrast, the formulations presented below are based upon 
  SOCLCP \eqref{eq.complementarity.form}. 
  \finbox
\end{remark}

A moment's consideration will show that \eqref{eq.complementarity.form} 
corresponds to the optimality condition of the following convex 
optimization problem: 
\begin{subequations}\label{P.truss.1}%
  \begin{alignat}{3}
    & \displaystyle 
    \text{Minimize}
    &{\quad}& \displaystyle
    \sum_{i=1}^{m} \Bigl( 
    q^{(t)}_{i} \Delta c_{\rr{e}_{i}}
    + \frac{1}{2} k_{i} \Delta c_{\rr{e}i}^{2}  \Bigr)
    + \sum_{i=1}^{m} \Bigl(
    R^{(t)}_{i} \Delta\gamma_{i} + \frac{1}{2} h_{i} \Delta \gamma_{i}^{2}  
    \Bigr)  
    - \bi{f}^{\top} \Delta\bi{u} 
    \label{P.truss.1.1} \\
    & \ST && \displaystyle
    \Delta c_{\rr{e}i} + \Delta c_{\rr{p}i} 
    = \bi{b}_{i}^{\top} \Delta\bi{u} , 
    \quad  i=1,\dots,m, 
    \label{P.truss.1.2} \\
    & && \displaystyle
    \Delta\gamma_{i} \ge | \Delta c_{\rr{p}i} | , 
    \quad  i=1,\dots,m . 
    \label{P.truss.1.3}
  \end{alignat}%
\end{subequations}%
It is worth noting that this problem is a variant of the total potential 
energy minimization formulation. 

Problem \eqref{P.truss.1} can be recast as (convex) QP and SOCP; 
see appendix~\ref{sec.SOCP} for reduction to SOCP. 
Therefore, it can be solved efficiently with a primal-dual 
interior-point method \citep{AL12,BtN01}. 
As an alternative approach, in this paper we examine a simpler 
gradient-based algorithm with acceleration. 

\begin{remark}
  The formulations presented above can readily be extended to the case 
  in which the tension and compression yield conditions are not 
  symmetric. 
  Suppose that the yield condition is given by 
  \begin{align*}
    \underline{q}_{i} \le q_{i} \le \overline{q}_{i} , 
  \end{align*}
  where $\underline{q}_{i}$ and $\overline{q}_{i}$ are constants. 
  This condition is equivalent to 
  \begin{align*}
    |q_{i} - \bar{\beta}_{i}| \le R_{i} 
  \end{align*}
  with $\bar{\beta}_{i}:=(\overline{q}_{i}+\underline{q}_{i})/2$ and 
  $R_{i}:=(\overline{q}_{i}-\underline{q}_{i})/2$. 
  Here, $\bar{\beta}_{i}$ and $R_{i}$ correspond to the center and 
  radius of the yield surface, respectively. 
  This modification from \eqref{eq.fundamental.yield.1} can be 
  realized by adding 
  \begin{align*}
    \sum_{i=1}^{m} \bar{\beta}_{i} \Delta c_{\rr{p}i} 
  \end{align*}
  to the objective function of problem \eqref{P.truss.1}. 
  A similar problem setting appears in section~\ref{sec.kinematic}. 
  \finbox
\end{remark}

\section{Accelerated proximal gradient method for elastoplastic analysis}
\label{sec.proximal}

In section~\ref{sec.unconstrained}, we reformulate the incremental problem 
as a form which is tractable within the framework of 
(accelerated) proximal gradient methods. 
A proximal gradient method and its accelerated version for this problem 
are presented in section~\ref{sec.not.accelerated} and 
section~\ref{sec.accelerated}, respectively. 

\subsection{Unconstrained formulation of minimum potential energy problem}
\label{sec.unconstrained}

In this section we recast problem \eqref{P.truss.1} as an unconstrained 
form. 

Since $R^{(t)}_{i}>0$ and $h_{i}>0$ $(i=1,\dots,m)$, 
all the inequality constraints of 
problem \eqref{P.truss.1} become active at the optimal solution. 
Therefore, by using $\Delta\gamma_{i} = | \Delta c_{\rr{p}i} |$ we can 
eliminate $\Delta\gamma_{i}$ as follows: 
\begin{subequations}\label{P.truss.2}%
  \begin{alignat}{3}
    & \displaystyle 
    \text{Minimize}
    &{\quad}& \displaystyle
    \sum_{i=1}^{m} \Bigl( 
    q^{(t)}_{i} \Delta c_{\rr{e}_{i}}
    + \frac{1}{2} k_{i} \Delta c_{\rr{e}i}^{2}  \Bigr)
    + \sum_{i=1}^{m} \Bigl(
    R^{(t)}_{i} |\Delta c_{\rr{p}i}| 
    + \frac{1}{2} h_{i} \Delta c_{\rr{p}i}^{2}  
    \Bigr)  
    - \bi{f}^{\top} \Delta\bi{u} \\
    & \ST && \displaystyle
    \Delta c_{\rr{e}i} + \Delta c_{\rr{p}i} 
    = \bi{b}_{i}^{\top} \Delta\bi{u} , 
    \quad  i=1,\dots,m . 
  \end{alignat}%
\end{subequations}%
Furthermore, by substituting the equality constraints to the objective 
function, we  can eliminate $\Delta c_{\rr{e}i}$ from \eqref{P.truss.2} 
as 
\begin{alignat}{2}
  &\Min
  &{\quad}&
  \sum_{i=1}^{m} \Bigl[
  q^{(t)}_{i} (\bi{b}_{i}^{\top} \Delta\bi{u} - \Delta c_{\rr{p}i})
  + \frac{1}{2} k_{i} (\bi{b}_{i}^{\top} \Delta\bi{u} - \Delta c_{\rr{p}i})^{2}  
  \Bigr] \notag\\
  &&{\quad}& \qquad
  {}+ \sum_{i=1}^{m} \Bigl(
  R^{(t)}_{i} |\Delta c_{\rr{p}i}| 
  + \frac{1}{2} h_{i} \Delta c_{\rr{p}i}^{2}  
  \Bigr)  
  - \bi{f}^{\top} \Delta\bi{u} . 
  \label{P.reduced.1}
\end{alignat}
This is an unconstrained convex optimization problem. 

In the following, for notational simplicity, we write problem 
\eqref{P.reduced.1} as 
\begin{alignat}{2}
  &\Min
  &{\quad}&
  \sum_{i=1}^{m} \Bigl[
  q^{(t)}_{i} (\bi{b}_{i}^{\top} \bi{v} - p_{i})
  + \frac{1}{2} k_{i} (\bi{b}_{i}^{\top} \bi{v} - p_{i})^{2}  
  \Bigr]\notag\\
  &&{\quad}& \qquad
  + \sum_{i=1}^{m} \Bigl(
  R^{(t)}_{i} |p_{i}| + \frac{1}{2} h_{i} p_{i}^{2}  
  \Bigr)  
  - \bi{f}^{\top} \bi{v}  
  \label{P.truss.3}
\end{alignat}
with $\bi{v}:=\Delta\bi{u}$ and $p_{i}:=\Delta c_{\rr{p}i}$ 
$(i=1,\dots,m)$. 
We propose to solve problem \eqref{P.truss.3} by applying 
an accelerated proximal gradient method. 

\begin{remark}
  Problem \eqref{P.truss.3} has a form very similar to the 
  $\ell_{1}$-regularized least-squares problem, 
  known as the LASSO \citep{Tib96,Tib11}. 
  The LASSO solves 
  \begin{align}
    \Min 
    {\quad}
    \| \ti{A} \bi{x} - \bi{b} \|^{2} 
    + \kappa \sum_{j=1}^{n} | x_{j} | , 
    \label{eq.LASSO.1}
  \end{align}
  where $\bi{x} \in \Re^{n}$ is a variable to be optimized and 
  $\kappa > 0$ is a constant regularization parameter. 
  A class of proximal gradient methods for solving problem 
  \eqref{eq.LASSO.1} is known as ISTA (iterative shrinkage-thresholding 
  algorithm); see, e.g., 
  \citet{CDLL98}, \citet{FN03}, \citet{DDM04}, and \citet{CW05}. 
  An accelerated version of ISTA due to \citet{BT09} is called 
  FISTA (fast iterative shrinkage-thresholding algorithm). 
  This paper is motivated by similarity between problems 
  \eqref{P.truss.3} and \eqref{eq.LASSO.1}; actually the algorithm 
  presented in section~\ref{sec.accelerated} is considered 
  essentially an application of FISTA to problem \eqref{P.truss.3}. 
  To the best of the author's knowledge, problem \eqref{P.truss.3} has 
  not been used in literature on computational mechanics. 
  \finbox
\end{remark}

\subsection{Proximal gradient method}
\label{sec.not.accelerated}

In this section we present a proximal gradient method for solving 
problem \eqref{P.truss.3}, which prepares the accelerated version 
appearing in section~\ref{sec.accelerated}. 

Define $g_{1} : \Re^{d} \times \Re^{m} \to \Re$ and 
$g_{2} : \Re^{m} \to \Re$ by 
\begin{align}
  g_{1}(\bi{v},\bi{p}) 
  &= \sum_{i=1}^{m} \Bigl[
  q^{(t)}_{i} (\bi{b}_{i}^{\top} \bi{v} - p_{i})
  + \frac{1}{2} k_{i} (\bi{b}_{i}^{\top} \bi{v} - p_{i})^{2}  
  \Bigr]
  + \sum_{i=1}^{m} \frac{1}{2} h_{i} p_{i}^{2}  
  - \bi{f}^{\top} \bi{v} , 
  \label{eq.proximal.setting.1} \\
  g_{2}(\bi{p}) 
  &= \sum_{i=1}^{m} R^{(t)}_{i} |p_{i}| , 
  \label{eq.proximal.setting.2}
\end{align}
which are closed proper convex functions. 
Particularly, $g_{1}$ is differentiable, 
and $\nabla g_{1}$ is Lipschitz continuous. 
We use $L$ to denote the Lipschitz constant of $\nabla g_{1}$. 
By making use of $g_{1}$ and $g_{2}$, problem \eqref{P.truss.3} can be 
written as 
\begin{align}
  \Min 
  {\quad} 
  g_{1}(\bi{v},\bi{p}) + g_{2}(\bi{p}) . 
  \label{P.proximal.decomposition}
\end{align}

A point $(\bi{v}^{*},\bi{p}^{*}) \in \Re^{d} \times \Re^{m}$ is optimal 
for problem \eqref{P.proximal.decomposition} if and only if 
\begin{align}
  \bi{0} &= \nabla_{\bi{v}}g_{1}(\bi{v}^{*},\bi{p}^{*}) , 
  \label{eq.optimality.proximal.4} \\
  \bi{0} &\in 
  \nabla_{\bi{p}}g_{1}(\bi{v}^{*},\bi{p}^{*}) 
  + \partial g_{2}(\bi{p}^{*}) , 
  \label{eq.optimality.proximal.3}
\end{align}
where 
\begin{align*}
  \nabla_{\bi{v}}g_{1} 
  = \pdif{g_{1}}{\bi{v}} , 
  \quad
  \nabla_{\bi{p}}g_{1}
  = \pdif{g_{1}}{\bi{p}} . 
\end{align*}
For any $\alpha > 0$, \eqref{eq.optimality.proximal.4} and 
\eqref{eq.optimality.proximal.3} hold if and only if the following 
equalities hold: 
\begin{align}
  \bi{v}^{*} 
  &= \bi{v}^{*} - \alpha \nabla_{\bi{v}} g_{1}(\bi{v}^{*},\bi{p}^{*}) , 
  \label{eq.optimality.proximal.5} \\
  \bi{p}^{*} 
  &= \bs{prox}_{\alpha g_{2}}
  (\bi{p}^{*} - \alpha\nabla_{\bi{p}} g_{1}(\bi{v}^{*},\bi{p}^{*})) . 
  \label{eq.optimality.proximal.6}
\end{align}
Equivalence of \eqref{eq.optimality.proximal.3} and 
\eqref{eq.optimality.proximal.6} follows from fundamental properties of 
the proximal mapping \citep{PB14}; see appendix~\ref{sec.fixed.point} 
for more accounts. 
The proximal gradient method can be designed from 
\eqref{eq.optimality.proximal.5} and \eqref{eq.optimality.proximal.6} as 
follows; see, e.g., \citet{CEg14} and \citet{PB14}. 

\begin{algorithm}\label{alg.truss.proximal.1}
  \hspace{1em}
  \begin{algstep}
    \setcounter{alnum}{-1}
    \item 
    Choose $\bi{v}_{0} \in \Re^{d}$, $\bi{p}_{0} \in \Re^{m}$, 
    $\alpha \in ]0, 1/L]$, and the 
    termination tolerance $\epsilon > 0$. 
    Set $l:=0$. 

    \item \label{alg.truss.proximal.1.step.1}
    Let 
    \begin{align*}
      \bi{v}_{l+1}
      &:= \bi{v}_{l} - 
      \alpha \nabla_{\bi{v}} g_{1}(\bi{v}_{l},\bi{p}_{l}) , \\
      \bi{p}_{l+1} 
      &:= \bs{prox}_{\alpha g_{2}}
      (\bi{p}_{l} - \alpha \nabla_{\bi{p}}g_{1}(\bi{v}_{l},\bi{p}_{l})) .
    \end{align*}
    
    \item 
    If 
    $\| (\bi{v}_{l}, \bi{p}_{l}) 
    - (\bi{v}_{l+1}, \bi{p}_{l+1}) \| \le \epsilon$, then terminate. 
    Otherwise, let $l \gets l+1$, and go to step~\ref{alg.truss.proximal.1.step.1}.
  \end{algstep}
\end{algorithm}

\subsection{Accelerated proximal gradient method with restart}
\label{sec.accelerated}

The convergence analysis of \refalg{alg.truss.proximal.1} can be found 
in literature, e.g., \citet{CEg14} and \citet{PB14}. 
It is known that $g_{1}(\bi{v}_{l},\bi{p}_{l})+g_{2}(\bi{p}_{l})$ 
converges to the optimal value at rate $1/l$. 
In this section we introduce the so-called FISTA-type improvement, due 
to \citet{BT09}, to achieve an accelerated convergence rate of 
order $1/l^{2}$. 
Also we incorporate the adaptive restart scheme proposed by \citet{OdC15} 
to ensure monotonic decrease of the objective function value. 

The following is an accelerated proximal gradient method 
with adaptive restart for the incremental elastoplastic analysis. 
\begin{algorithm}\label{alg.truss.proximal.3}
  \hspace{1em}
  \begin{algstep}
    \setcounter{alnum}{-1}
    \item \label{alg.truss.proximal.3.step0}
    Choose $\bi{v}_{0} \in \Re^{d}$, $\bi{p}_{0} \in \Re^{m}$, 
    $\alpha \in ]0,1/L]$, and 
    the termination tolerance $\epsilon > 0$. 
    Set $l:=1$, $\bi{\mu}_{1} := \bi{v}_{0}$, 
    $\bi{\rho}_{1} := \bi{p}_{0}$, and $\tau_{1}:=1$. 

    \item \label{alg.truss.proximal.3.step2}
    Let 
    \begin{align*}
      \bi{v}_{l}
      &:= \bi{\mu}_{l} - \alpha 
      \nabla_{\bi{\mu}}g_{1}(\bi{\mu}_{l}, \bi{\rho}_{l}) , \\
      \bi{p}_{l}
      &:= \bs{prox}_{\alpha g_{2}} 
      (\bi{\rho}_{l} - \alpha 
      \nabla_{\bi{\rho}}g_{1}(\bi{\mu}_{l}, \bi{\rho}_{l})) . 
    \end{align*}
    
    \item \label{alg.truss.proximal.3.step3}
    Let 
    \begin{align*}
      \tau_{l+1}  
      := \frac{1}{2} \Bigl( 1 + \sqrt{1 + 4 \tau_{l}^{2}} \Bigr) .
    \end{align*}
    
    \item \label{alg.truss.proximal.3.step4}
    If $g_{1}(\bi{v}_{l},\bi{p}_{l}) + g_{2}(\bi{p}_{l}) 
      < g_{1}(\bi{v}_{l-1},\bi{p}_{l-1}) + g_{2}(\bi{p}_{l-1})$, then let 
    \begin{align*}
      \bi{\mu}_{l+1}  
      &:= \bi{v}_{l} 
      + \frac{\tau_{l}-1}{\tau_{l+1}} (\bi{v}_{l} - \bi{v}_{l-1}) , \\
      \bi{\rho}_{l+1}  
      &:= \bi{p}_{l} + \frac{\tau_{l}-1}{\tau_{l+1}} 
      (\bi{p}_{l} - \bi{p}_{l-1}) . 
    \end{align*}
    Otherwise, let 
    \begin{align*}
      \tau_{l+1}  &:=  1 , \\
      \bi{\mu}_{l+1}  &:= \bi{v}_{l} , \\
      \bi{\rho}_{l+1}  &:= \bi{p}_{l} . 
    \end{align*}
    
    \item
    If $\| (\bi{v}_{l}, \bi{p}_{l}) 
    - (\bi{v}_{l-1}, \bi{p}_{l-1}) \| \le \epsilon$, then terminate. 
    Otherwise, let $l \gets l+1$, and go to step~\ref{alg.truss.proximal.3.step2}.
  \end{algstep}
\end{algorithm}

\begin{remark}
  If we replace step~\ref{alg.truss.proximal.3.step3} with 
  $\tau_{l+1}:=1$, then \refalg{alg.truss.proximal.3} reverts to 
  \refalg{alg.truss.proximal.1}. 
  \finbox
\end{remark}

Computation of step~\ref{alg.truss.proximal.3.step2} can be carried out 
in an efficient manner as follows. 

We begin with computation of the proximal mapping of $\alpha g_{2}$ with 
$\alpha > 0$, which is defined as  
\begin{align}
  \bs{prox}_{\alpha g_{2}}(\bi{w}) 
  = \argmin_{\bi{z}} \Bigl\{
  \alpha \sum_{i=1}^{m} R^{(t)}_{i} |z_{i}| 
  + \frac{1}{2} \| \bi{z} - \bi{w} \|^{2}
  \Bigr\} . 
  \label{eq.ell1.proximal.map.1}
\end{align}
Since we have that 
\begin{align}
  \min_{\bi{z}} \Bigl\{
  \alpha \sum_{i=1}^{m} R^{(t)}_{i} |z_{i}| 
  + \frac{1}{2} \| \bi{z} - \bi{w} \|^{2}
  \Bigr\} 
  = \sum_{i=1}^{m} \min_{z_{i}} \Bigl\{
  \alpha R^{(t)}_{i} |z_{i}| + \frac{1}{2} (z_{i} - w_{i})^{2}
  \Bigr\} , 
  \label{eq.proximal.computation.3}
\end{align}
the optimal $\bi{z}$ in the right-hand side of 
\eqref{eq.ell1.proximal.map.1} can be found by solving the set of the 
one-dimensional optimization problems in the right-hand side of 
\eqref{eq.proximal.computation.3}. 
Then, it is known that the optimal solution for each $i$ can be obtained 
via the soft-threshold function (also known as the shrinkage operator) as 
\begin{align}
  \argmin_{z_{i}} \Bigl\{
  \alpha R^{(t)}_{i} |z_{i}| + \frac{1}{2} (z_{i} - w_{i})^{2}
  \Bigr\} = 
  \begin{dcases*}
    0 
    & if $|w_{i}| \le \alpha R^{(t)}_{i}$, \\
    w_{i} - \alpha R^{(t)}_{i} \sign(w_{i}) 
    & otherwise; 
  \end{dcases*}
  \label{eq.proximal.computation.1}
\end{align}
see, e.g., \citet{CW05}, \citet{BT09}, \citet{CEg14}, and \citet{PB14}. 
Consequently, we obtain 
\begin{align}
  \bs{prox}_{\alpha g_{2}}(\bi{w}) = 
  \diag(\sign(\bi{w})) \max \{ |\bi{w}| - \alpha \bi{R}^{(t)}, \bi{0} \} .
  \label{eq.prox.operator.of.g2}
\end{align}

We next consider computation of $\nabla g_{1}$. 
Define $\ti{B} \in \Re^{m \times d}$ by 
\begin{align}
  \ti{B} = 
  \begin{bmatrix}
    \bi{b}_{1}^{\top} \\
    \vdots \\
    \bi{b}_{m}^{\top} \\
  \end{bmatrix}
  , 
\end{align}
which is the compatibility matrix. 
By using this notation, definition \eqref{eq.proximal.setting.1} of 
$g_{1}$ yields 
\begin{align}
  \nabla_{\bi{v}} g_{1}(\bi{v},\bi{p}) 
  &= \ti{B}^{\top} \bi{q}^{(t)} 
  + \ti{B}^{\top} \diag(\bi{k}) \ti{B} \bi{v} 
  - \ti{B}^{\top} \diag(\bi{k}) \bi{p} - \bi{f}  , 
  \label{eq.gradient.g1.of.v} \\
  \nabla_{\bi{p}} g_{1}(\bi{v},\bi{p}) 
  &= -\bi{q}^{(t)}_{i} 
  + \diag(\bi{k}) (\bi{p} - \ti{B} \bi{v}) 
  + \diag(\bi{h}) \bi{p} . 
  \label{eq.gradient.g1.of.p}
\end{align}
For computing \eqref{eq.gradient.g1.of.v} and \eqref{eq.gradient.g1.of.p}, 
it is convenient to define $\bi{e} \in \Re^{m}$ by 
\begin{align}
  \bi{e} 
  = \ti{B} \bi{v} - \bi{p} . 
  \label{eq.recover.c.e}
\end{align}
It is worth noting that $\bi{e}$ corresponds to the incremental elastic 
elongation, $\Delta\bi{c}_{\rr{e}}$, in problem \eqref{P.truss.1}. 
By using $\bi{e}$, \eqref{eq.gradient.g1.of.v} and 
\eqref{eq.gradient.g1.of.p} can be calculated as 
\begin{align}
  \nabla_{\bi{v}} g_{1}(\bi{v},\bi{p}) 
  &= \ti{B}^{\top} (\diag(\bi{k}) \bi{e} + \bi{q}^{(t)}) 
  - \bi{f} ,  
  \label{eq.gradient.f.sub.u} \\
  \nabla_{\bi{p}} g_{1}(\bi{v},\bi{p}) 
  &= \diag(\bi{h}) \bi{p}  - \diag(\bi{k}) \bi{e} - \bi{q}^{(t)} . 
  \label{eq.gradient.f.sub.p}
\end{align}
Consequently, by using 
\eqref{eq.prox.operator.of.g2}, \eqref{eq.recover.c.e}, 
\eqref{eq.gradient.f.sub.u}, and \eqref{eq.gradient.f.sub.p}, we see 
that computation at step~\ref{alg.truss.proximal.3.step2} of 
\refalg{alg.truss.proximal.3} can be performed as follows: 
\begin{align}
  \bi{\varepsilon}_{l} 
  &:= \ti{B} \bi{\mu}_{l} - \bi{\rho}_{l} , 
  \label{eq.accelerated.update.1} \\
  \bi{v}_{l} 
  &:= \bi{\mu}_{l} - \alpha 
  [ \ti{B}^{\top} (\diag(\bi{k}) \bi{\varepsilon}_{l} + \bi{q}^{(t)}) 
  - \bi{f} ] , 
  \label{eq.accelerated.update.2} \\
  \bi{w}_{l}
  &:= \bi{\rho}_{l} - \alpha 
  (\diag(\bi{h}) \bi{\rho}_{l} - \diag(\bi{k}) \bi{\varepsilon}_{l} 
  - \bi{q}^{(t)}) , 
  \label{eq.accelerated.update.3} \\
  \bi{p}_{l} 
  &:= \diag(\sign(\bi{w}_{l})) \max \{ |\bi{w}_{l}| - \alpha \bi{R}^{(t)}, \bi{0} \}
  . 
  \label{eq.accelerated.update.4}
\end{align}
Here, $\bi{\varepsilon}_{l}$ and $\bi{w}_{l}$ are auxiliary variables. 

\begin{remark}
  The most expensive part of \refalg{alg.truss.proximal.3} is 
  computation at step~\ref{alg.truss.proximal.3.step2}. 
  As seen in \eqref{eq.accelerated.update.1}, 
  \eqref{eq.accelerated.update.2}, \eqref{eq.accelerated.update.3}, and 
  \eqref{eq.accelerated.update.4}, 
  this essentially amounts to two matrix-vector products and 
  four component-wise vector products. 
  Here, $\ti{B}$ is a sparse matrix, and hence the two matrix-vector 
  products may exploit this sparsity effectively. 
  \finbox
\end{remark}

\begin{remark}
  \refalg{alg.truss.proximal.3} does not contain any process of solving a 
  system of linear equations. 
  Therefore, \refalg{alg.truss.proximal.3} does not require any 
  linear-equations solver. 
  If a conventional method for elastoplastic analysis is applied to 
  large-scale problems, then a linear-equations solver dominates the 
  computational cost. 
  Hence, usually iterative methods are used for solving the equilibrium 
  equation with the tangent stiffness matrix. 
  Also, parallel computing, 
  such as domain decomposition methods \citep{NAD07,CKSV14}, is often required. 
  When we adopt an approach based upon mathematical programming, an 
  interior-point method solves a system of linear equations to find the 
  search direction at each iteration. 
  Hence, to solve a large-scale problem an iterative solver is usually 
  employed for computation of the search direction; 
  see, e.g., \citet{JSS00}, \citet{PRVJ00}, \citet{KKLBF07}, and \citet{BGVZ07}. 
  In contrast, \refalg{alg.truss.proximal.3} does not use a 
  linear-equations solver at all. 
  In other words, \refalg{alg.truss.proximal.3} is explicit (and also 
  simple), although it solves a problem discretized with 
  a fully implicit Euler scheme. 
  \finbox
\end{remark}

\begin{remark}
  Like other approaches based upon mathematical programming, 
  \refalg{alg.truss.proximal.3} does not resort to a consistent tangent 
  stiffness matrix. 
  Moreover, it does not require any procedure to determine whether each 
  member undergoes plastic loading or elastic unloading. 
  \finbox
\end{remark}

At step~\ref{alg.truss.proximal.3.step0} of 
\refalg{alg.truss.proximal.3}, we can determine the step size, $\alpha$, 
as follows. 
From \eqref{eq.gradient.g1.of.v} and \eqref{eq.gradient.g1.of.p}, the 
Hessian matrix of $g_{1}$ can be obtained as 
\begin{align}
  \nabla^{2} g_{1}(\bi{v},\bi{p}) 
  &= 
  \left[\begin{array}{@{}c|c@{\,}}
   \displaystyle
     \ti{B}^{\top} \diag(\bi{k}) \ti{B} 
     & {-}\ti{B}^{\top} \diag(\bi{k}) \\
     \hline
     {-}\diag(\bi{k}) \ti{B} 
     & \diag(\bi{k} + \bi{h}) \\
  \end{array}
  \right]  
  \label{eq.Hessian.g1.1} \\
  &= 
  \left[\begin{array}{@{}c|c@{\,}}
   \ti{B}^{\top} & \ti{O} \\
   \hline
     -\ti{I} & \ti{I} \\
  \end{array}
  \right] 
  \left[\begin{array}{@{}c|c@{\,}}
   \diag(\bi{k}) & \ti{O} \\
   \hline
     \ti{O} & \diag(\bi{h}) \\
  \end{array}
  \right] 
  \left[\begin{array}{@{}c|c@{\,}}
   \ti{B} & -\ti{I} \\
   \hline
     \ti{O} & \ti{I} \\
  \end{array}
  \right] . 
  \label{eq.Hessian.g1.2}
\end{align}
Recall that $k_{i}>0$ and $h_{i}>0$ $(i=1,\dots,m)$. 
Moreover, for a stable (more precisely, kinematically determinate) truss, 
$\ti{B}$ is of row full rank. 
Therefore, from \eqref{eq.Hessian.g1.2} we see that 
$\nabla^{2} g_{1}(\bi{v},\bi{p})$ is positive definite, which implies 
that $g_{1}$ is strongly convex. 
Furthermore, the maximum eigenvalue of 
$\nabla^{2} g_{1}(\bi{v},\bi{p})$ is equal to $L$, 
i.e., the Lipschitz constant of $\nabla g_{1}$. 
One obvious choice for determining $\alpha$ is, therefore, to find the 
maximum eigenvalue of the matrix in \eqref{eq.Hessian.g1.1} and set 
$\alpha:=1/L$. 
Another choice is to find an upper bound for $L$ that can be performed 
much faster than its exact value. 
For notational convenience, let 
$H = (H_{ij}) = \nabla^{2}g_{1}(\bi{v},\bi{p})$. 
It follows from the Gershgorin disc theorem that $L'$ defined by 
\begin{align}
  L' = \max \Bigl\{
  H_{ii} + \sum_{j \not= i} |H_{ij}| 
  \Bigm|
  i=1,\dots,d+m
  \Bigr\}
  \label{eq.def.upper.bound}
\end{align}
satisfies $L' \ge L$; see, e.g., \citet[Theorem~6.1.1]{HJ85}. 
Then we may set $\alpha:=1/L'$.

\section{Mixed isotropic/kinematic hardening}
\label{sec.kinematic}

In this section we consider a plasticity model that combines linear 
isotropic hardening and linear kinematic hardening. 

We begin by formulating the incremental problem. 
To incorporate the kinematical hardening, the yield condition, 
\eqref{eq.fundamental.yield.1}, is replaced with 
\begin{align*}
  |q_{i} - \beta_{i}|  \le R_{i} . 
\end{align*}
Here, $\beta_{i} \in \Re$ is an internal force corresponding to the back 
stress. 
Let $\theta \in [0,1]$ be a constant. 
The evolutions of $R_{i}$ and $\beta_{i}$ are given by 
\begin{align}
  R_{i} 
  &=  R^{(t)}_{i} + \theta h_{i} \Delta\gamma_{i} , 
  \label{eq.kinematic.evolution.1} \\
  \beta_{i} 
  &= \beta^{(t)}_{i} + (1-\theta) h_{i} \Delta c_{\rr{p}i} . 
  \label{eq.kinematic.evolution.2}
\end{align}
Here, $\theta$ is the ratio of the effect of isotropic hardening to the 
total strain hardening. 
Particularly, $\theta=1$ corresponds to the pure isotropic hardening, and 
$\theta=0$ corresponds to the pure kinematic hardening. 
Thus, \eqref{eq.fundamental.truss.7} in section~\ref{sec.fundamental} is 
replaced with \eqref{eq.kinematic.evolution.1}, and 
\eqref{eq.kinematic.evolution.2} is newly added. 
Consequently, the incremental problem can be formulated as 
\begin{subequations}\label{eq.complementarity.kinematic}
  \begin{alignat}{2}
    & \Delta c_{\rr{e}i} + \Delta c_{\rr{p}i} 
    = \bi{b}_{i}^{\top} \Delta\bi{u} , 
    &{\quad}& i=1,\dots,m, \\
    & q_{i} = q^{(t)}_{i} + k_{i} \Delta c_{\rr{e}i} , 
    &{\quad}& i=1,\dots,m, \\
    & \beta_{i} = \beta^{(t)}_{i} + (1-\theta) h_{i} \Delta c_{\rr{p}i} , 
    &{\quad}& i=1,\dots,m, \\
    & \sum_{i=1}^{m} q_{i} \bi{b}_{i} = \bi{f} , \\
    & R^{(t)}_{i} + \theta h_{i} \Delta\gamma_{i} \ge |q_{i} - \beta_{i}| , 
    \quad
    \Delta\gamma_{i} \ge |\Delta c_{\rr{p}i}| ,  \notag\\
    & \qquad
    \begin{bmatrix}
      R^{(t)}_{i} + \theta h_{i} \Delta\gamma_{i} \\ q_{i} - \beta_{i} \\
    \end{bmatrix}
    ^{\top}
    \begin{bmatrix}
      \Delta\gamma_{i} \\ {-}\Delta c_{\rr{p}i} \\
    \end{bmatrix}
    = 0 , 
    &{\quad}& i=1,\dots,m. 
  \end{alignat}
\end{subequations}
Like \eqref{eq.complementarity.form} in section~\ref{sec.fundamental}, 
this is an SOCLCP.

It is easy to verify that \eqref{eq.complementarity.kinematic} is 
the optimality condition of the following convex optimization problem: 
\begin{subequations}\label{P.kinematic.1}
  \begin{alignat}{3}
    & \displaystyle 
    \text{Minimize}
    &{\quad}& \displaystyle
    \sum_{i=1}^{m} \Bigl( 
    q^{(t)}_{i} \Delta c_{\rr{e}_{i}}
    + \frac{1}{2} k_{i} \Delta c_{\rr{e}i}^{2}  \Bigr) \notag\\
    &&{\quad}& \displaystyle
    \qquad
    + \sum_{i=1}^{m} \Bigl[
    R^{(t)}_{i} \Delta\gamma_{i} 
    + \frac{1}{2} \theta h_{i} \Delta \gamma_{i}^{2}  
    + \beta^{(t)}_{i} \Delta c_{\rr{p}i} 
    + \frac{1}{2} (1-\theta) h_{i} \Delta c_{\rr{p}i}^{2}  
    \Bigr]  \notag\\
    &&{\quad}& \displaystyle
    \qquad
    - \bi{f}^{\top} \Delta\bi{u} 
    \label{P.kinematic.1.1} \\
    & \ST && \displaystyle
    \Delta c_{\rr{e}i} + \Delta c_{\rr{p}i} 
    = \bi{b}_{i}^{\top} \Delta\bi{u} , 
    \quad  i=1,\dots,m, \\
    & && \displaystyle
    \Delta\gamma_{i} \ge | \Delta c_{\rr{p}i} | , 
    \quad  i=1,\dots,m. 
    \label{P.kinematic.1.3}
  \end{alignat}%
\end{subequations}%
Since $R^{(t)}_{i} > 0$ and $\theta h_{i} \ge 0$ $(i=1,\dots,m)$, the 
constraints in \eqref{P.kinematic.1.3} become active at the optimal 
solution. 
Therefore, without changing the optimal solution we can substitute 
$\Delta c_{\rr{p}i}^{2} = \Delta \gamma_{i}^{2}$ into 
\eqref{P.kinematic.1.1}. 
This results in 
\begin{align}
  \sum_{i=1}^{m} \Bigl( 
  q^{(t)}_{i} \Delta c_{\rr{e}_{i}}
  + \frac{1}{2} k_{i} \Delta c_{\rr{e}i}^{2}  \Bigr) 
  + \sum_{i=1}^{m} \Bigl(
  R^{(t)}_{i} \Delta\gamma_{i} 
  + \frac{1}{2} h_{i} \Delta \gamma_{i}^{2}  
  + \beta^{(t)}_{i} \Delta c_{\rr{p}i} 
  \Bigr)  
  - \bi{f}^{\top} \Delta\bi{u} , 
  \label{P.kinematic.3}
\end{align}
which is similar to \eqref{P.truss.1.1}. 
More precisely, the difference is only the 
presence of $\beta^{(t)}_{i} \Delta c_{\rr{p}i}$ $(i=1,\dots,m)$. 
Therefore, \refalg{alg.truss.proximal.3} can be applied in a very 
similar manner. 
Namely, we just replace \eqref{eq.accelerated.update.3} used 
at step~\ref{alg.truss.proximal.3.step2} with 
\begin{align*}
  \bi{w}_{l}
  &:= \bi{\rho}_{l} - \alpha 
  (\diag(\bi{h}) \bi{\rho}_{l} + \bi{\beta}^{(t)} - \diag(\bi{k}) \bi{\varepsilon}_{l} 
  - \bi{q}^{(t)})  
\end{align*}
to find the solution. 
Subsequently, $R^{(t)}_{i}$ and $\beta^{(t)}_{i}$ should be updated 
according to \eqref{eq.kinematic.evolution.1} and 
\eqref{eq.kinematic.evolution.2} for computation of the next time increment.

\section{Piecewise-linear hardening}
\label{sec.piecewise}

In this section we consider an accelerated proximal gradient method for 
a piecewise-linear hardening model. 
It is worth noting that QP formulations for piecewise-linear model is 
known in literature, e.g., \citet{Mai68}. 

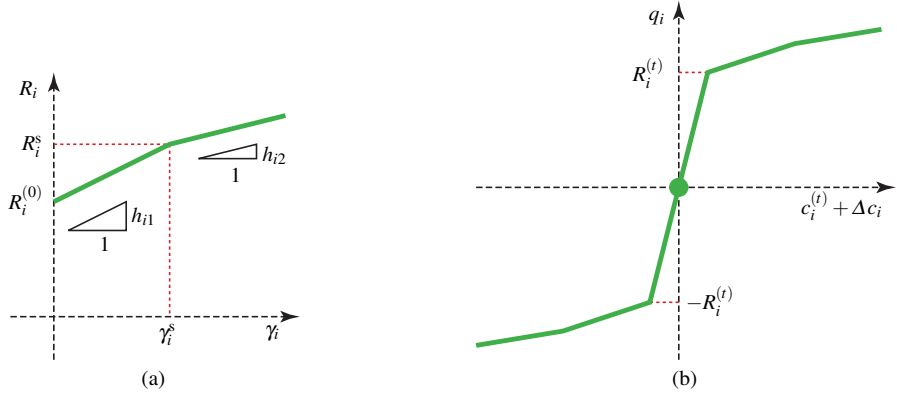
\begin{figure}[tp]
  \centering
  \subfigure[]{
  \label{fig.plastic_law_4}
  \includegraphics[scale=0.60]{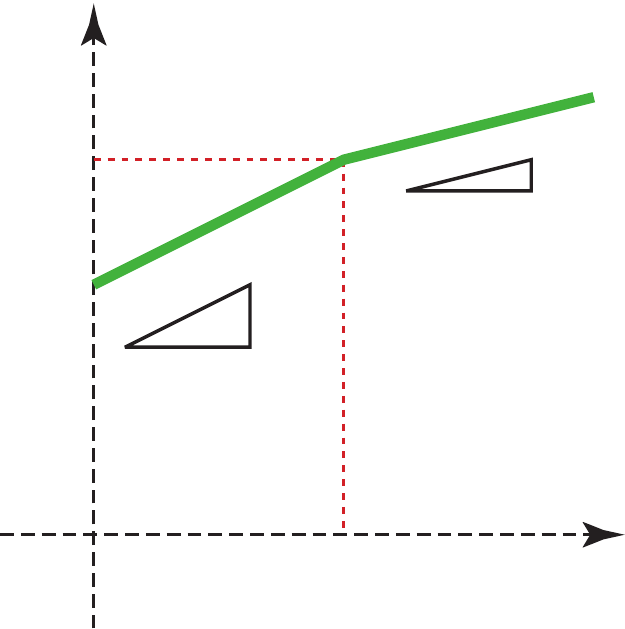}
  \begin{picture}(0,0)
    \put(-180,0){
    \put(165,9){{\small $\gamma_{i}$}}
    \put(73,100){{\small $R_{i}$}}
    \put(69,57){{\small $R^{(0)}_{i}$}}
    \put(74,79){{\small $R^{\rr{s}}_{i}$}}
    \put(125,8){{\small $\gamma^{\rr{s}}_{i}$}}
    \put(103,41){{\small $1$}}
    \put(115,52){{\small $h_{i1}$}}
    \put(152,67){{\small $1$}}
    \put(164,75){{\small $h_{i2}$}}
    }
  \end{picture}
  }
  \hfill
  \subfigure[]{
  \label{fig.plastic_law_3}
  \centering
  \includegraphics[scale=0.60]{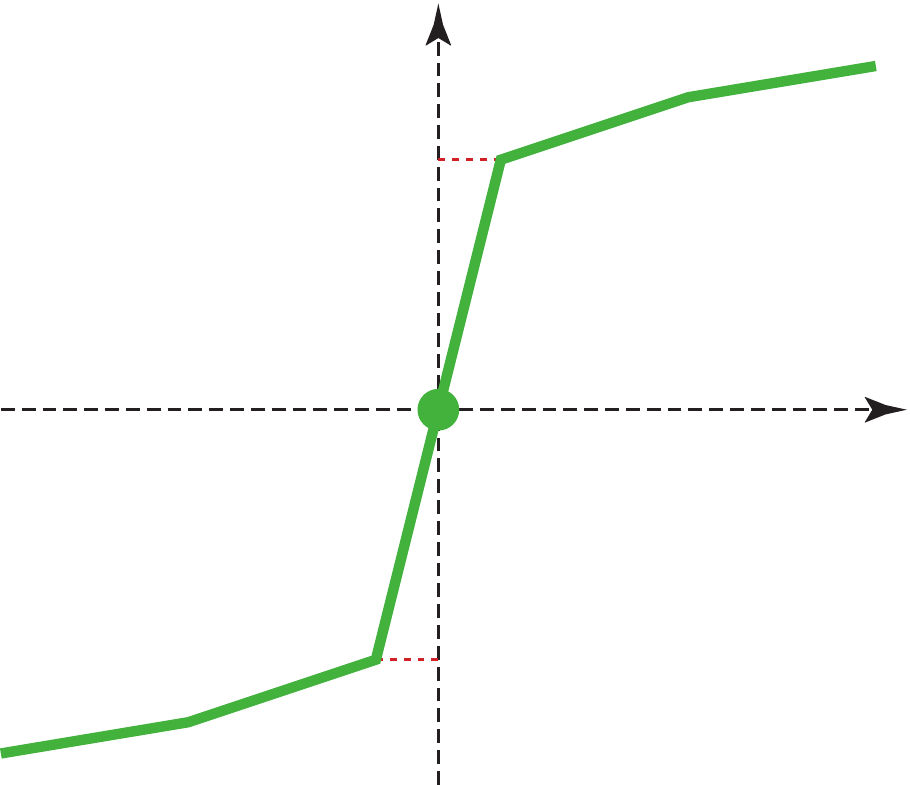}
  \begin{picture}(0,0)
    \put(-180,-10){
    \put(85,138){{\small $q_{i}$}}
    \put(142,65){{\small $c^{(t)}_{i} + \Delta c_{i}$}}
    \put(99,29){{\small $-R^{(t)}_{i}$}}
    \put(79,115){{\small $R^{(t)}_{i}$}}
    }
  \end{picture}
  }
  \caption{A piecewise-linear hardening model. 
  \subref{fig.plastic_law_3} The evolution of $R_{i}$; and 
  \subref{fig.plastic_law_3} the axial force versus elongation relation.}
  \label{fig.plastic_law_piecewise}
\end{figure}

Suppose that the evolution of $R_{i}$, with respect to the plastic 
multiplier, is given as illustrated in 
\reffig{fig.plastic_law_piecewise}. 
That is, when $R_{i}$ attains at $R^{\rr{s}}_{i}$, which is a given 
positive constant, the hardening modulus decreases from $h_{i1}$ to 
$h_{i2}$, where $h_{i1} > 0$ and $h_{i2} \in ]0,h_{i1}[$ 
are given constants. 
Under this hypothesis we formulate a quasi-static incremental problem. 

Recall that, in section~\ref{sec.fundamental}, we have decomposed the 
incremental elongation by \eqref{eq.fundamental.truss.1}. 
Instead, in this section we consider the following decomposition: 
\begin{align}
  \Delta c_{i} 
  = \Delta c_{\rr{e}i} + \Delta c_{\rr{p}i} + \Delta c_{\rr{s}i} . 
\end{align}
Here, $\Delta c_{\rr{s}i}$ is a variable used to assess the plastic 
elongation after the axial force attains $R^{\rr{s}}_{i}$; see 
\eqref{eq.piecewise.complementarity.1} and 
\eqref{eq.piecewise.complementarity.3} for more precise interpretation. 
Define $\Delta\gamma_{i1}$ and $\Delta\gamma_{i2}$ by 
\begin{align}
  \Delta\gamma_{i1}  &= |\Delta c_{\rr{p}i}| , \\
  \Delta\gamma_{i2}  &= |\Delta c_{\rr{s}i}| . 
  \label{eq.piecewise.complementarity.1}
\end{align}
The evolution of $R_{i}$ defined as \reffig{fig.plastic_law_piecewise} 
can be written in terms of $\Delta\gamma_{i1}$ and $\Delta\gamma_{i2}$ 
as 
\begin{alignat}{3}
  & R_{i} \le R^{\rr{s}}_{i} 
  &{\quad}& \Rightarrow \quad
  R_{i} = R^{(t)}_{i} + h_{i1} \Delta\gamma_{i1} , \ 
  \Delta\gamma_{i2} = 0 , 
  \label{eq.piecewise.complementarity.4} \\
  & R_{i} > R^{\rr{s}}_{i} 
  &{\quad}& \Rightarrow \quad
  R_{i} = R^{\rr{s}}_{i} + h_{i2} (\Delta\gamma_{i1} + \Delta\gamma_{i2}) .
  \label{eq.piecewise.complementarity.5}
\end{alignat}
For simplicity, define $\eta_{i}$ by 
\begin{align*}
  \eta_{i} = \frac{h_{i1}h_{i2}}{h_{i1}-h_{i2}} , 
\end{align*}
which is a positive constant. 
A moment's consideration will show that 
\eqref{eq.piecewise.complementarity.4} and 
\eqref{eq.piecewise.complementarity.5} are equivalent to 
\begin{align}
  R_{i}
  &= R^{(t)}_{i} + h_{i1} \Delta\gamma_{i1} + \eta_{i} \Delta\gamma_{i2} 
\end{align}
and 
\begin{alignat}{2}
  & R_{i} \le R^{\rr{s}}_{i} 
  &{\quad}& \Rightarrow \quad
  \Delta\gamma_{i2} = 0 , 
  \label{eq.piecewise.complementarity.2} \\
  & R_{i} > R^{\rr{s}}_{i} 
  &{\quad}& \Rightarrow \quad
  R_{i} = R^{\rr{s}}_{i} + \eta_{i} \Delta\gamma_{i2} . 
  \label{eq.piecewise.complementarity.3}
\end{alignat}
It can be readily verified that 
\eqref{eq.piecewise.complementarity.1}, 
\eqref{eq.piecewise.complementarity.2}, and 
\eqref{eq.piecewise.complementarity.3} are equivalent to the following 
second-order cone complementarity condition. 
\begin{proposition}\label{prop.piecewise}
  Assume $R_{i} > 0$. 
  Then $R_{i}$, $\Delta\gamma_{i2}$, and $\Delta c_{\rr{s}i}$ satisfy 
  \begin{alignat*}{3}
    & R_{i} \le R^{\rr{s}}_{i} 
    &{\quad}& \Rightarrow \quad 
    \Delta\gamma_{i2} = \Delta c_{\rr{s}i} = 0 , \\
    & R_{i} > R^{\rr{s}}_{i} 
    &{\quad}& \Rightarrow \quad 
    \Delta\gamma_{i2} = |\Delta c_{\rr{s}i}| , 
    \quad
    R_{i} = R^{\rr{s}}_{i} + \eta_{i} \Delta\gamma_{i2} 
  \end{alignat*}
  if and only if they satisfy 
  \begin{align*}
    R^{\rr{s}}_{i} + \eta_{i} \Delta\gamma_{i2} \ge |R_{i}| , 
    \quad
    \Delta\gamma_{i2} \ge |\Delta c_{\rr{s}i}| ,  
    \quad
    \begin{bmatrix}
      R^{\rr{s}}_{i} + \eta_{i} \Delta\gamma_{i2} \\ R_{i} \\
    \end{bmatrix}
    ^{\top}
    \begin{bmatrix}
      \Delta\gamma_{i2} \\ -\Delta c_{\rr{s}i} \\
    \end{bmatrix}
    = 0 . 
  \end{align*}
\end{proposition}

\refprop{prop.complementarity} and \refprop{prop.piecewise} combine to 
give the following formulation of the incremental problem: 
\begin{subequations}\label{eq.complementarity.piecewise}
\begin{alignat}{2}
  & \Delta c_{\rr{e}i} + \Delta c_{\rr{p}i} + \Delta c_{\rr{s}i} 
  = \bi{b}_{i}^{\top} \Delta\bi{u} , 
  &{\quad}& i=1,\dots,m, \\
  & q_{i} = q^{(t)}_{i} + k_{i} \Delta c_{\rr{e}i} , 
  &{\quad}& i=1,\dots,m, \\
  & \sum_{i=1}^{m} q_{i} \bi{b}_{i} = \bi{f} , \\
  & R^{(t)}_{i} + h_{i1} \Delta\gamma_{i1} \ge |q_{i}| , 
  \quad
  \Delta\gamma_{i1} \ge |\Delta c_{\rr{p}i}| , \notag\\
  &{\qquad} 
  \begin{bmatrix}
    R^{(t)}_{i} + h_{i1} \Delta\gamma_{i1} \\ q_{i} \\
  \end{bmatrix}
  ^{\top}
  \begin{bmatrix}
    \Delta\gamma_{i1} \\ -\Delta c_{\rr{p}i} \\
  \end{bmatrix}
  = 0 , 
  &{\quad}& i=1,\dots,m, \\
  & R^{\rr{s}}_{i} + \eta_{i} \Delta\gamma_{i2} \ge |q_{i}| , 
  \quad
  \Delta\gamma_{i2} \ge |\Delta c_{\rr{s}i}| ,  \notag\\
  &{\qquad} 
  \begin{bmatrix}
    R^{\rr{s}}_{i} + \eta_{i} \Delta\gamma_{i2} \\ q_{i} \\
  \end{bmatrix}
  ^{\top}
  \begin{bmatrix}
    \Delta\gamma_{i2} \\ -\Delta c_{\rr{s}i} \\
  \end{bmatrix}
  = 0 , 
  &{\quad}& i=1,\dots,m. 
\end{alignat}
\end{subequations}
We can show that \eqref{eq.complementarity.piecewise} corresponds to the 
optimality condition of the following convex optimization problem: 
\begin{subequations}\label{eq.P.piecewise.2}
  \begin{alignat}{3}
    & \displaystyle 
    \text{Minimize}
    &{\quad}& \displaystyle
    \sum_{i=1}^{m} \Bigl( 
    q^{(t)}_{i} \Delta c_{\rr{e}_{i}}
    + \frac{1}{2} k_{i} \Delta c_{\rr{e}i}^{2}  \Bigr)  \notag\\
    &&{\quad}& \displaystyle {\qquad}
    + \sum_{i=1}^{m} \Bigl(
    R^{(t)}_{i} \Delta\gamma_{i1} 
    + \frac{1}{2} h_{i1} \Delta \gamma_{i1}^{2}  
    + R^{\rr{s}}_{i} \Delta\gamma_{i2} 
    + \frac{1}{2} \eta_{i} \Delta \gamma_{i2}^{2}  
    \Bigr)  \notag\\
    &&{\quad}& \displaystyle {\qquad}
    - \bi{f}^{\top} \Delta\bi{u} \\
    & \ST && \displaystyle
    \Delta c_{\rr{e}i} + \Delta c_{\rr{p}i} + \Delta c_{\rr{s}i} 
    = \bi{b}_{i}^{\top} \Delta\bi{u} , 
    \quad  i=1,\dots,m, \\
    & && \displaystyle
    \Delta\gamma_{i1} \ge | \Delta c_{\rr{p}i} | , 
    \quad  i=1,\dots,m, \\
    & && \displaystyle
    \Delta\gamma_{i2} \ge | \Delta c_{\rr{s}i} | , 
    \quad  i=1,\dots,m. 
  \end{alignat}%
\end{subequations}%

In a manner similar to section~\ref{sec.unconstrained}, we can recast 
problem \eqref{eq.P.piecewise.2} as an unconstrained nonsmooth convex 
optimization problem. 
Firstly, since $R^{(t)}_{i} > 0$, $R^{\rr{s}}_{i} > 0$, $h_{i1} > 0$, 
and $\eta_{i} > 0$ $(i=1,\dots,m)$, 
all the inequality constraints of problem \eqref{eq.P.piecewise.2} 
become active at an optimal solution. 
Therefore, 
$\Delta\gamma_{i1}$ and $\Delta\gamma_{i2}$ can be eliminated. 
Next, by making use of the equality constraints, 
we can eliminate $\Delta c_{\rr{e}i}$. 
As a result, we see that problem \eqref{eq.P.piecewise.2} is equivalent 
to 
\begin{align}
  \Min 
  {\quad} & 
  \sum_{i=1}^{m} \Bigl( 
  q^{(t)}_{i} (\bi{b}_{i}^{\top} \Delta\bi{u} - \Delta c_{\rr{p}i} - \Delta c_{\rr{s}i})
  + \frac{1}{2} k_{i} 
  (\bi{b}_{i}^{\top} \Delta\bi{u} - \Delta c_{\rr{p}i} - \Delta c_{\rr{s}i})^{2}  
  \Bigr)  \notag\\
  &{\quad} 
  + \sum_{i=1}^{m} \Bigl(
  R^{(t)}_{i} |\Delta c_{\rr{p}i}| 
  + \frac{1}{2} h_{i1} \Delta c_{\rr{p}i}^{2}  
  + R^{\rr{s}}_{i} |\Delta c_{\rr{s}i}| 
  + \frac{1}{2} \eta_{i} \Delta c_{\rr{s}i}^{2}  
  \Bigr)  
  - \bi{f}^{\top} \Delta\bi{u} . 
  \label{eq.P.piecewise.Delta}
\end{align}
For notational simplicity, we write problem \eqref{eq.P.piecewise.Delta} 
as 
\begin{align}
  \Min 
  {\quad} & 
  \sum_{i=1}^{m} \Bigl( 
  q^{(t)}_{i} (\bi{b}_{i}^{\top} \bi{v} - p_{i} - s_{i})
  + \frac{1}{2} k_{i} 
  (\bi{b}_{i}^{\top} \bi{v} - p_{i} - s_{i})^{2}  
  \Bigr)  \notag\\
  &{\quad} 
  + \sum_{i=1}^{m} \Bigl(
  R^{(t)}_{i} |p_{i}| + \frac{1}{2} h_{i1} p_{i}^{2}  
  + R^{\rr{s}}_{i} |s_{i}| + \frac{1}{2} \eta_{i} s_{i}^{2}  
  \Bigr)  
  - \bi{f}^{\top} \bi{v} 
  \label{eq.P.piecewise.1}
\end{align}
with $\bi{v}:= \Delta\bi{u}$, $p_{i}:=\Delta c_{\rr{p}i}$, and 
$s_{i}:=\Delta c_{\rr{s}i}$ $(i=1,\dots,m)$. 

Define convex functions 
$g_{1}: \Re^{d} \times \Re^{m} \times \Re^{m} \to \Re$ 
and $g_{2}: \Re^{m} \times \Re^{m} \to \Re$ by 
\begin{align}
  g_{1}(\bi{v},\bi{p},\bi{s}) 
  &= \sum_{i=1}^{m} \Bigl[
  q^{(t)}_{i} (\bi{b}_{i}^{\top} \bi{v} - p_{i} - s_{i})
  + \frac{1}{2} k_{i} 
  (\bi{b}_{i}^{\top} \bi{v} - p_{i} - s_{i})^{2}  
  \Bigr]  \notag\\
  & \qquad 
  + \sum_{i=1}^{m} \frac{1}{2} h_{i1} p_{i}^{2} 
  + \sum_{i=1}^{m} \frac{1}{2} \eta_{i} s_{i}^{2} 
  - \bi{f}^{\top} \bi{v} , 
  \label{eq.P.piecewise.g.1} \\
  g_{2}(\bi{p},\bi{s}) 
  &= \sum_{i=1}^{m} R^{(t)}_{i} |p_{i}| 
  + \sum_{i=1}^{m} R^{\rr{s}}_{i} |s_{i}| . 
  \label{eq.P.piecewise.g.2}
\end{align}
Here, $g_{1}$ is differentiable, and $\nabla g_{1}$ is 
Lipschitz continuous. 
Problem \eqref{eq.P.piecewise.1} can be written as 
\begin{align}
  \Min
  {\quad}
  g_{1}(\bi{v},\bi{p},\bi{s}) + g_{2}(\bi{p},\bi{s}) . 
  \label{eq.P.piecewise.abstract.1}
\end{align}
Then we can design an accelerated proximal gradient method in a manner 
similar to section~\ref{sec.accelerated}; details appear in 
appendix~\ref{sec.algorithm_piece}.

\section{Numerical experiments}
\label{sec.ex}
The presented algorithms were implemented with MATLAB ver.\ 8.4.0. 
Comparison is performed with \textsf{QUADPROG} \citep{MATLAB}, 
\textsf{IPOPT} ver.\ 3.11.3 \citep{WA06} via 
the Matlab interface \citep{IPOPT.for.MATLAB}, and 
\textsf{PATH} ver.\ 4.7.03 \citep{DF95,FM00}
via the Matlab Interface \citep{FM99}. 
\textsf{QUADPROG} is a MATLAB built-in function for convex 
quadratic programming (QP). 
We use an implementation of an interior-point method 
by setting the parameter \texttt{Algorithm} to \texttt{interior-point-convex}. 
\textsf{IPOPT} is a primal-dual interior-point method with a filter 
line-search method for nonlinear programming. 
We set the parameters 
\texttt{options.ipopt.hessian\_approximation} and 
\texttt{options.ipopt.tol} to 
\texttt{limited-memory} and $10^{-1}$, respectively. 
\textsf{PATH} is a nonsmooth Newton method to solve mixed 
complementarity problems. 
We apply \textsf{PATH} to solve the KKT condition for the QP problem in 
\eqref{P.truss.1}. 
The termination criterion of \refalg{alg.truss.proximal.3} and 
\refalg{alg.piecewise} is 
$\| (\bi{v}_{l}, \bi{p}_{l}) 
 - (\bi{v}_{l-1}, \bi{p}_{l-1}) \|_{\infty} \le \epsilon$ 
with $\epsilon=10^{-8}\,\mathrm{m}$. 
Computation was carried out on 
a $2.6\,\mathrm{GHz}$ Intel Core i5 processor with 
$8\,\mathrm{GB}$ RAM. 

\begin{figure}[tp]
  \centering
  \subfigure[]{
  \label{fig.model_10_perspective}
  \includegraphics[scale=0.50]{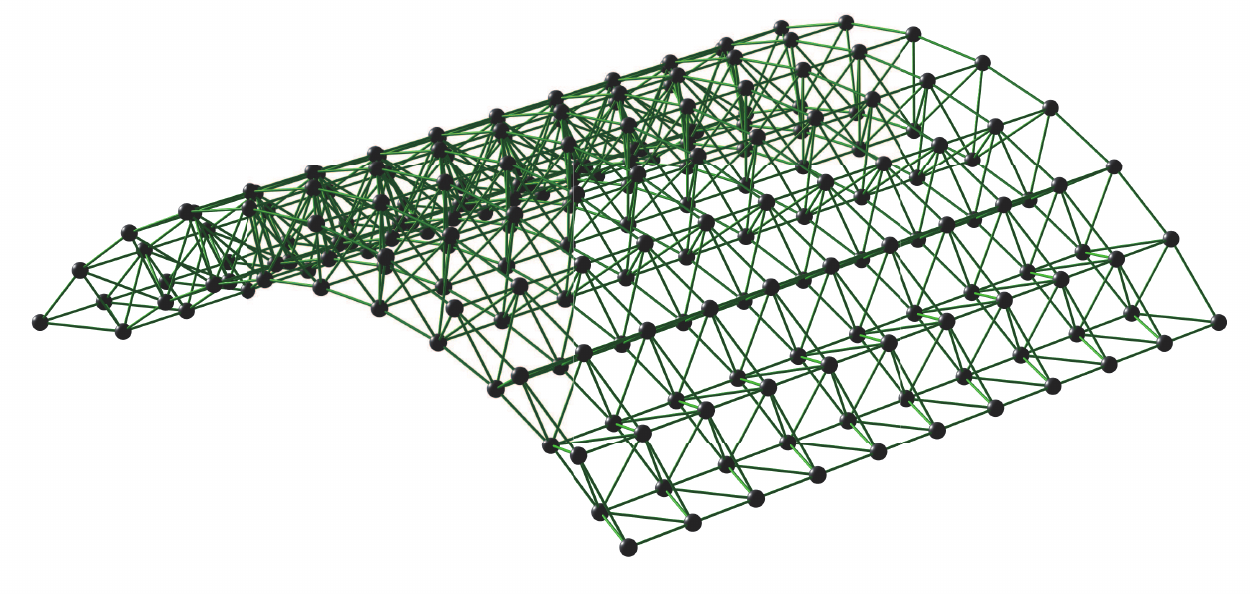}
  }
  \hfill
  \subfigure[]{
  \label{fig.model_10_z_axis}
  \scalebox{0.80}{
  \includegraphics[scale=0.5]{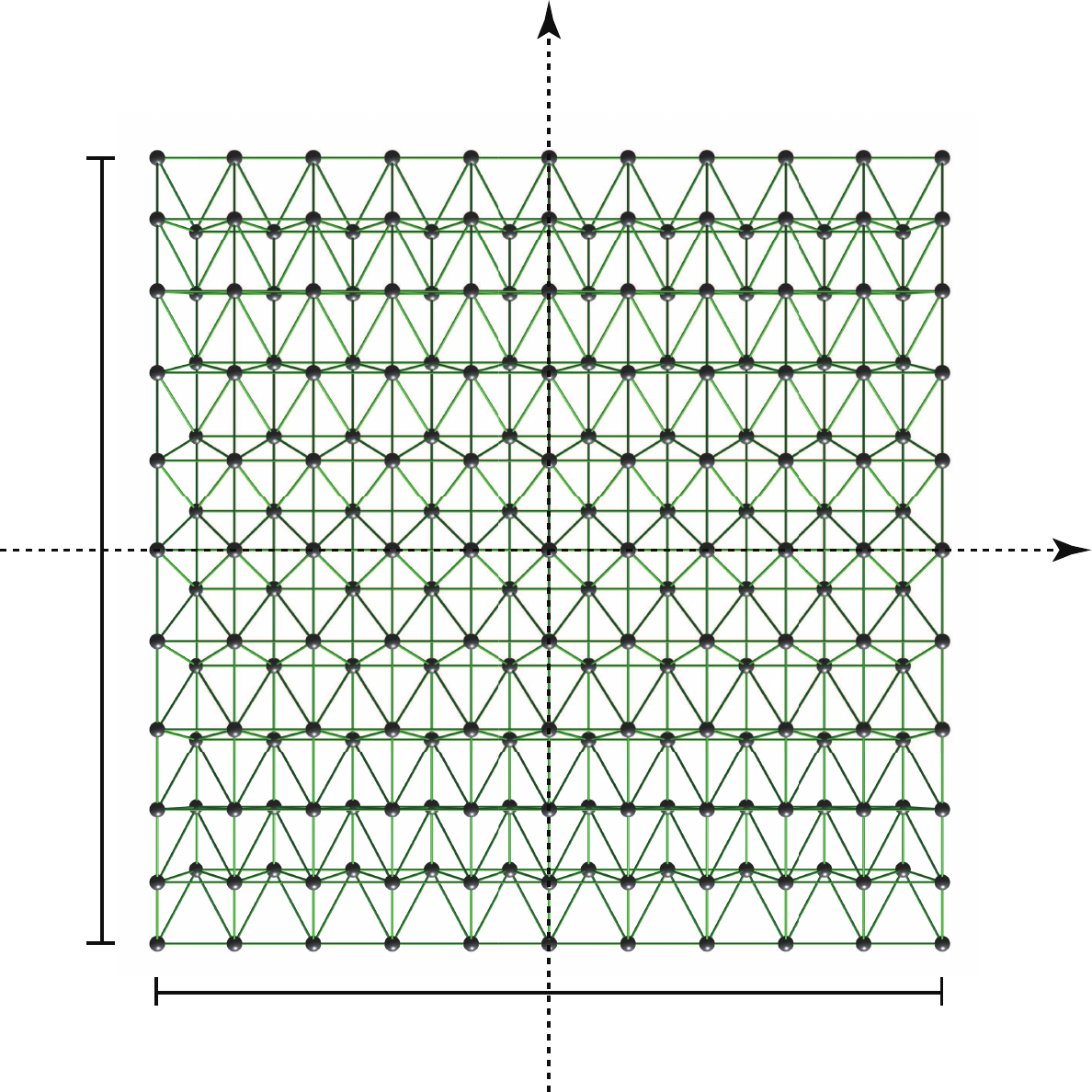}
  \begin{picture}(0,0)
    \put(-180,-10){
    \put(105,17){{\small $N_{X}${\,}@$1\,\mathrm{m}$}}
    \put(-13,120){{\small $N_{Y}${\,}@$1\,\mathrm{m}$}}
    \put(165,85){{\small $X$}}
    \put(83,173){{\small $Y$}}
    }
  \end{picture}
  }
  }
  \par\bigskip\bigskip
  \subfigure[]{
  \label{fig.model_10_y_axis}
  \scalebox{0.80}{
  \includegraphics[scale=0.65]{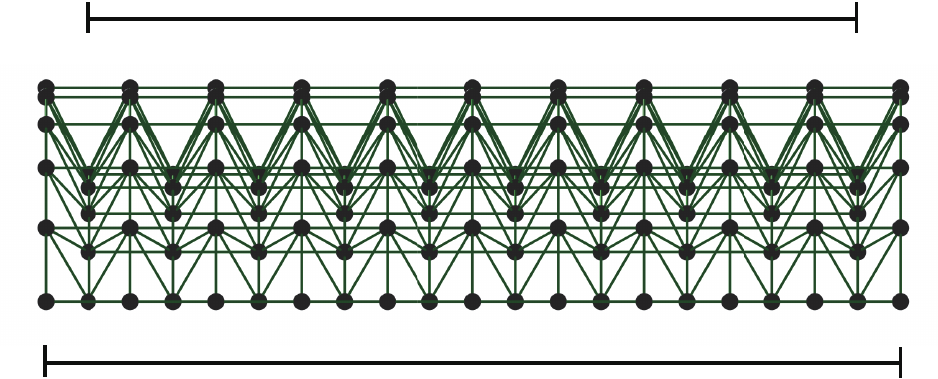}
  \begin{picture}(0,0)
    \put(-260,-40){
    \put(182,34){{\small $N_{X}${\,}@$1\,\mathrm{m}$}}
    \put(182,112){{\small $(N_{X}-1)${\,}@$1\,\mathrm{m}$}}
    }
  \end{picture}
  }
  }
  \hfill
  \subfigure[]{
  \label{fig.model_10_x_axis}
  \scalebox{0.80}{
  \includegraphics[scale=0.65]{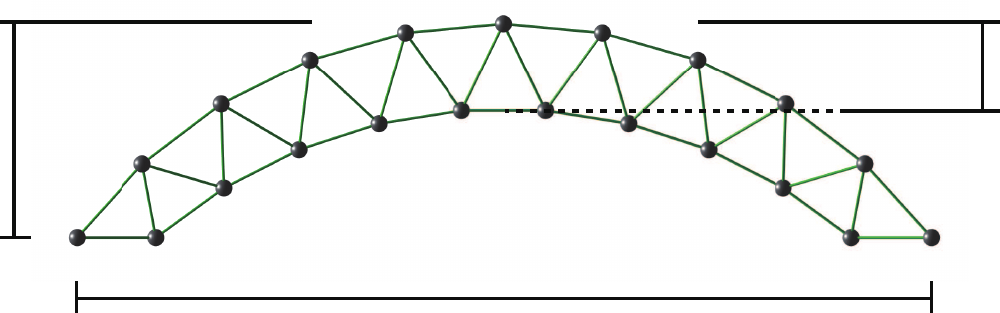}
  \begin{picture}(0,0)
    \put(-280,-40){
    \put(165,46){{\small $N_{Y}${\,}@$1\,\mathrm{m}$}}
    \put(278,83){{\small $1\,\mathrm{m}$}}
    \put(62,70){{\small $N_{Y}/4\,\mathrm{m}$}}
    }
  \end{picture}
  }
  }
  \caption{The problem setting with $(N_{X},N_{Y})=(10,10)$.
  \subref{fig.model_10_perspective} Perspective view; 
  \subref{fig.model_10_z_axis} plan; 
  \subref{fig.model_10_y_axis} elevation: and 
  \subref{fig.model_10_x_axis} side view. 
  }
  \label{fig.model_10}
\end{figure}

In the following numerical experiments, we consider a truss shown in 
\reffig{fig.model_10}, where $N_{X}$ and $N_{Y}$ are varied to generate 
problem instances with diverse sizes. 
This barrel vault truss is a two-way space grid with square pyramids. 
In the direction of the $X$-axis, the nodes are aligned with 
regular intervals, as shown in \reffig{fig.model_10_y_axis}. 
In the $Y Z$-plane, as seen in \reffig{fig.model_10_x_axis}, the top 
layer nodes are aligned on a circle equiangularly. 
Also, the bottom layer nodes are equiangularly aligned on a circle with 
the same radius as the one of the top layer nodes. 
All the lowest nodes of the top layer are pin-supported. 
The number of members, $m$, and the number of degrees of freedom of 
displacements, $d$, are listed in \reftab{tab.instances}. 
The cross-sectional area of each member is 
$a_{i}=500\,\mathrm{mm^{2}}$ and Young's modulus is 
$E=200\,\mathrm{GPa}$. 


\begin{table}[bp]
  \centering
  \caption{Characteristics of the problem instances.}
  \label{tab.instances}
  \begin{tabular}{lrr}
    \toprule
    $(N_{X},N_{Y})$ & $m$ & $d$  \\
    \midrule
    $(10,10)$ & 800 & 597 \\
    $(20,20)$ & 3{,}200  & 2{,}397 \\
    $(30,30)$ & 7{,}200  & 5{,}397 \\
    $(40,40)$ & 12{,}800 & 9{,}597 \\
    $(50,50)$ & 20{,}000 & 14{,}997 \\
    $(60,60)$ & 28{,}800 & 21{,}597 \\
    $(70,70)$ & 39{,}200 & 29{,}397 \\
    $(80,80)$ & 51{,}200 & 38{,}397 \\
    $(90,90)$ & 64{,}800 & 48{,}597 \\
    $(100,100)$ & 80{,}000  & 59{,}997 \\
    $(110,110)$ & 96{,}800  & 72{,}597 \\
    $(120,120)$ & 115{,}200 & 86{,}397 \\
    $(130,130)$ & 135{,}200 & 101{,}397 \\
    \bottomrule
  \end{tabular}
\end{table}

\begin{remark}
  It is well known that vectorizing MATLAB code often increases 
  computational efficiency drastically \citep{MATLAB}. 
  All the computations at step~\ref{alg.truss.proximal.3.step2} of 
  \refalg{alg.truss.proximal.3}, i.e., 
  \eqref{eq.accelerated.update.1}, \eqref{eq.accelerated.update.2}, 
  \eqref{eq.accelerated.update.3}, and \eqref{eq.accelerated.update.4}, 
  can be implemented in vectorized forms. 
  Namely, a component-wise vector product can be carried out with the  
  MATLAB function \texttt{times}. 
  Also, for calculation of \eqref{eq.accelerated.update.4}, we can apply 
  \texttt{abs}, \texttt{max}, and \texttt{sign} functions to vectors. 
  Similarly, step~\ref{alg.piecewise.step.1} of 
  \refalg{alg.piecewise} can also be implemented in 
  vectorized forms. 
  \finbox
\end{remark}

\subsection{Holonomic (path-independent) analysis}

In sections~\ref{sec.ex.1} and~\ref{sec.ex.2} we perform holonomic 
analysis, i.e., we assume that no elastic unloading takes place and 
consider a relatively large loading step. 
With these examples we attempt to evaluate the efficiency of the 
presented algorithm when it is applied to moderately large-scale problem 
instances. 
Two variants of the step size discussed in section~\ref{sec.accelerated} 
are examined. 
In the following, by \textsf{APGM-E} we mean that the step size is set 
to $\alpha:=1/L$ with Lipschitz constant $L$ of $\nabla g_{1}$, and 
by \textsf{APGM-G} we mean that $\alpha:=1/L'$ with $L'$ defined by 
\eqref{eq.def.upper.bound}. 

We set the current axial forces to 
$q^{(0)}_{i}=0$ $(i=1,\dots,m)$. 
The hardening moduli are $h_{i}=0.1 k_{i}$ $(i=1,\dots,m)$. 
The yielding stress is $\sigma_{\rr{y}}=200\,\mathrm{MPa}$ and, 
accordingly, 
$R^{(0)}_{i}=\sigma_{\rr{y}} a_{i} = 100\,\mathrm{kN}$ $(i=1,\dots,m)$. 

\subsubsection{Example (I): linear hardening model}
\label{sec.ex.1}

\begin{table}[bp]
  \centering
  \caption{The computational results of the proposed methods in 
  example (I).}
  \label{tab.result.ex.1.1}
  \begin{tabular}{lrrrrrrr}
    \toprule 
    & \multicolumn{4}{c}{\textsf{APGM-E}} 
    & \multicolumn{3}{c}{\textsf{APGM-G}} \\
    \cmidrule(r){2-5} \cmidrule(l){6-8} 
    $(N_{X},N_{Y})$ 
    & Iter.\ & Time (s) & \texttt{eigs} (s) & Rel.\ diff.\
    & Iter.\ & Time (s) & Rel.\ diff.\ \\
    \midrule
    $(10,10)$ 
    & 404 & 0.1 & (0.0) & $-8.6$e-12 
    & 350 & 0.1 & 1.2e-08 \\
    $(20,20)$ 
    & 549 & 0.3 & (0.1) & 5.3e-08 
    & 732 & 0.3 & 5.0e-08  \\
    $(30,30)$ 
    & 669 & 0.8 & (0.3) & 6.1e-07 
    & 893 & 0.6 & 6.7e-07 \\
    $(40,40)$ 
    & 1,450 & 2.4 & (0.9) & 7.9e-09 
    & 1,927 & 2.1 & 8.9e-09 \\
    $(50,50)$ 
    & 2,112 & 5.2 & (1.7) & 2.4e-08 
    & 2,778 & 4.5 & 2.5e-08 \\
    $(60,60)$ 
    & 3,700 & 10.1 & (1.5) & 1.6e-07 
    & 4,900 & 11.2 & 1.6e-07 \\
    $(70,70)$ 
    & 7,772 & 27.5 & (4.1) & $-1.7$e-09 
    & 4,853 & 14.6 & 2.3e-05 \\
    $(80,80)$ 
    & 6,237 & 30.0 & (6.3) & 5.2e-06 
    & 7,075 & 27.0 & 3.9e-05 \\
    $(90,90)$ 
    & 8,784 & 53.0 & (11.8) & 1.9e-07 
    & 9,235 & 43.3 & 9.5e-06 \\
    $(100,100)$ 
    & 11,158 & 84.7 & (17.9) & 1.4e-08 
    & 14,724 & 88.2 & 1.6e-08 \\
    $(110,110)$ 
    & 10,807 & 101.8 & (20.8) & 1.5e-05 
    & 14,251 & 107.1 & 1.5e-05 \\
    $(120,120)$ 
    & 12,973 & 159.8 & (41.5) & 1.4e-05 
    & 17,097 & 156.3 & 1.4e-05 \\
    $(130,130)$ 
    & 15,309 & 253.5 & (45.3) & 8.2e-06 
    & 20,165 & 269.0 & 8.2e-06 \\
    \bottomrule
  \end{tabular}
\end{table}

\begin{table}[bp]
  \centering
  \caption{The computational results of \textsf{QUADPROG}, 
  \textsf{IPOPT}, and \textsf{PATH} in example (I)}
  \label{tab.result.ex.1.2}
  \begin{tabular}{lrrrrrrrr}
    \toprule 
    & \multicolumn{2}{c}{\textsf{QUADPROG}} 
    & \multicolumn{3}{c}{\textsf{IPOPT}} 
    & \multicolumn{3}{c}{\textsf{PATH}} \\
    \cmidrule(r){2-3} \cmidrule(lr){4-6} \cmidrule(l){7-9} 
    $(N_{X},N_{Y})$ 
    & Iter.\ & Time (s) 
    & Iter.\ & Time (s) & Rel.\ diff.\
    & Iter.\ & Time (s) & Rel.\ diff.\ \\
    \midrule
    $(10,10)$ 
    & 10 & 0.4 
    & 36 & 0.7 & 4.0e-06 
    & 12 & 0.5 & $-4.8$e-10 \\
    $(20,20)$ 
    & 11 & 1.7 
    & 49 & 3.8 & 3.6e-05 
    & 12 & 8.7 & $-3.1$e-10 \\
    $(30,30)$ 
    & 11 & 6.2 
    & 74 & 14.8 & 2.9e-05 
    & 12 & 95.0 & $-8.6$e-10 \\
    $(40,40)$ 
    & 11 & 12.8 
    & 192 & 70.3 & 4.5e-06 
    & 12 & 699.2 & $-6.2$e-09 \\
    $(50,50)$ 
    & 11 & 26.2 
    & 423 & 236.4 & 1.8e-06 
    & --- & ($>1200.0$) & --- \\
    $(60,60)$ 
    & 11 & 44.9 
    & 508 & 425.5 & 1.5e-05 
    & --- & ($>1200.0$) & --- \\
    $(70,70)$ 
    & 11 & 74.0 
    & 958 & 1,116.7 & 1.3e-04 
    & --- & ($>1200.0$) & --- \\
    $(80,80)$ 
    & 11 & 112.2 
    & --- & ($>1200.0$) & --- 
    & --- & ($>1200.0$) & --- \\
    $(90,90)$ 
    & 11 & 167.0 
    & --- & ($>1200.0$) & --- 
    & --- & ($>1200.0$) & --- \\
    $(100,100)$ 
    & 11 & 241.3 
    & --- & ($>1200.0$) & --- 
    & --- & ($>1200.0$) & --- \\
    $(110,110)$ 
    & 11 & 353.4 
    & --- & ($>1200.0$) & --- 
    & --- & ($>1200.0$) & --- \\
    $(120,120)$ 
    & 11 & 495.8 
    & --- & ($>1200.0$) & --- 
    & --- & ($>1200.0$) & --- \\
    $(130,130)$ 
    & 11 & 633.3 
    & --- & ($>1200.0$) & --- 
    & --- & ($>1200.0$) & --- \\
    \bottomrule
  \end{tabular}
\end{table}

\begin{figure}[tp]
  \centering
  \includegraphics[scale=0.40]{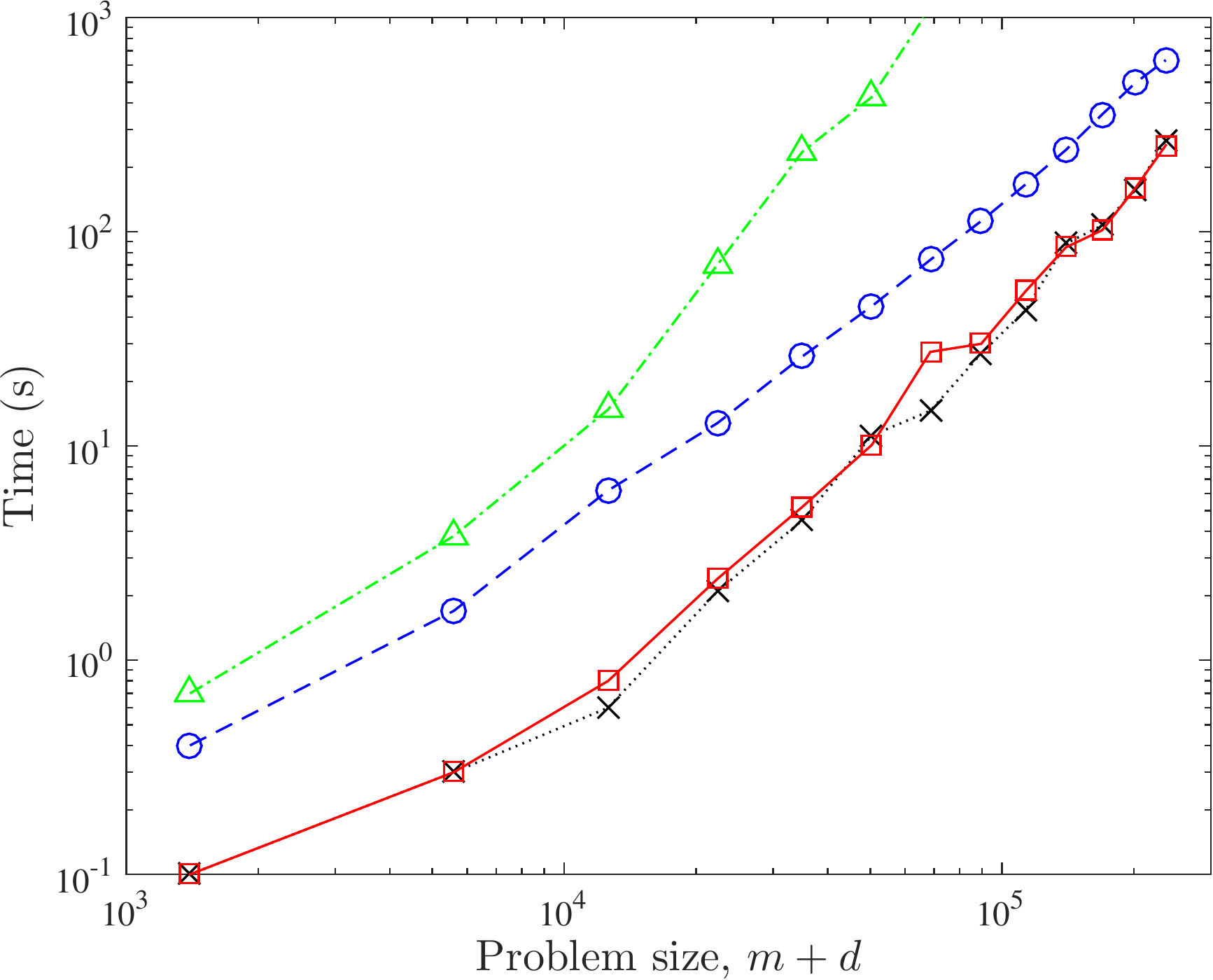}
  \caption{The computational time of example (I). 
  ``{\footnotesize $\square$}'' \textsf{APGM-E}; 
  ``{\footnotesize $\times$}''  \textsf{APGM-G}; 
  ``{\footnotesize $\circ$}'' \textsf{QUADPROG}; and 
  ``{\footnotesize $\vartriangle$}'' \textsf{IPOPT}. 
  }
  \label{fig.holonomic_time}
\end{figure}

\begin{figure}[tp]
  \centering
  \includegraphics[scale=0.40]{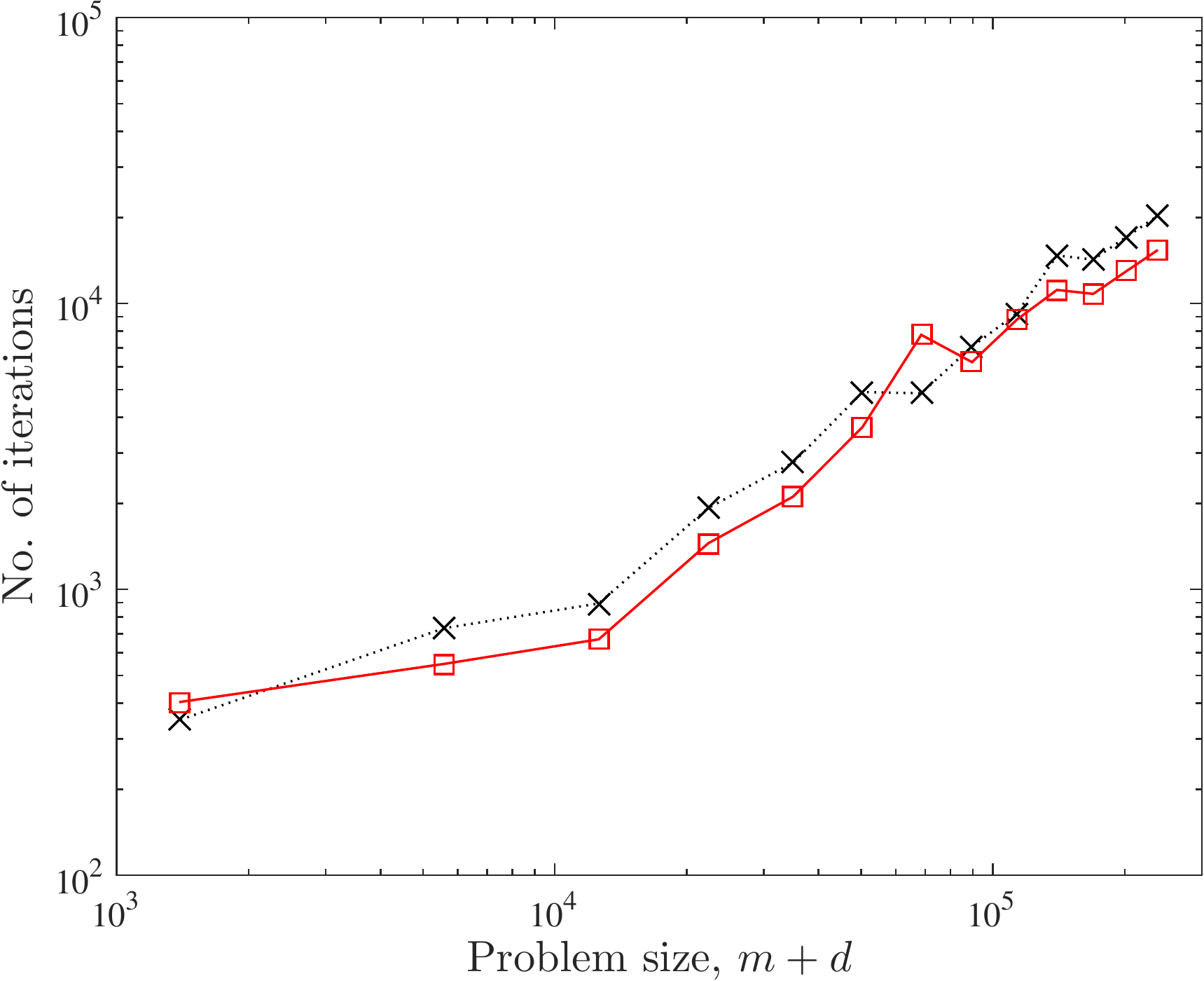}
  \caption{The number of iterations of example (I). 
  ``{\footnotesize $\square$}'' \textsf{APGM-E}; and 
  ``{\footnotesize $\times$}''  \textsf{APGM-G}. 
  }
  \label{fig.holonomic_iter}
\end{figure}

\begin{figure}[tp]
  \centering
  \includegraphics[scale=0.40]{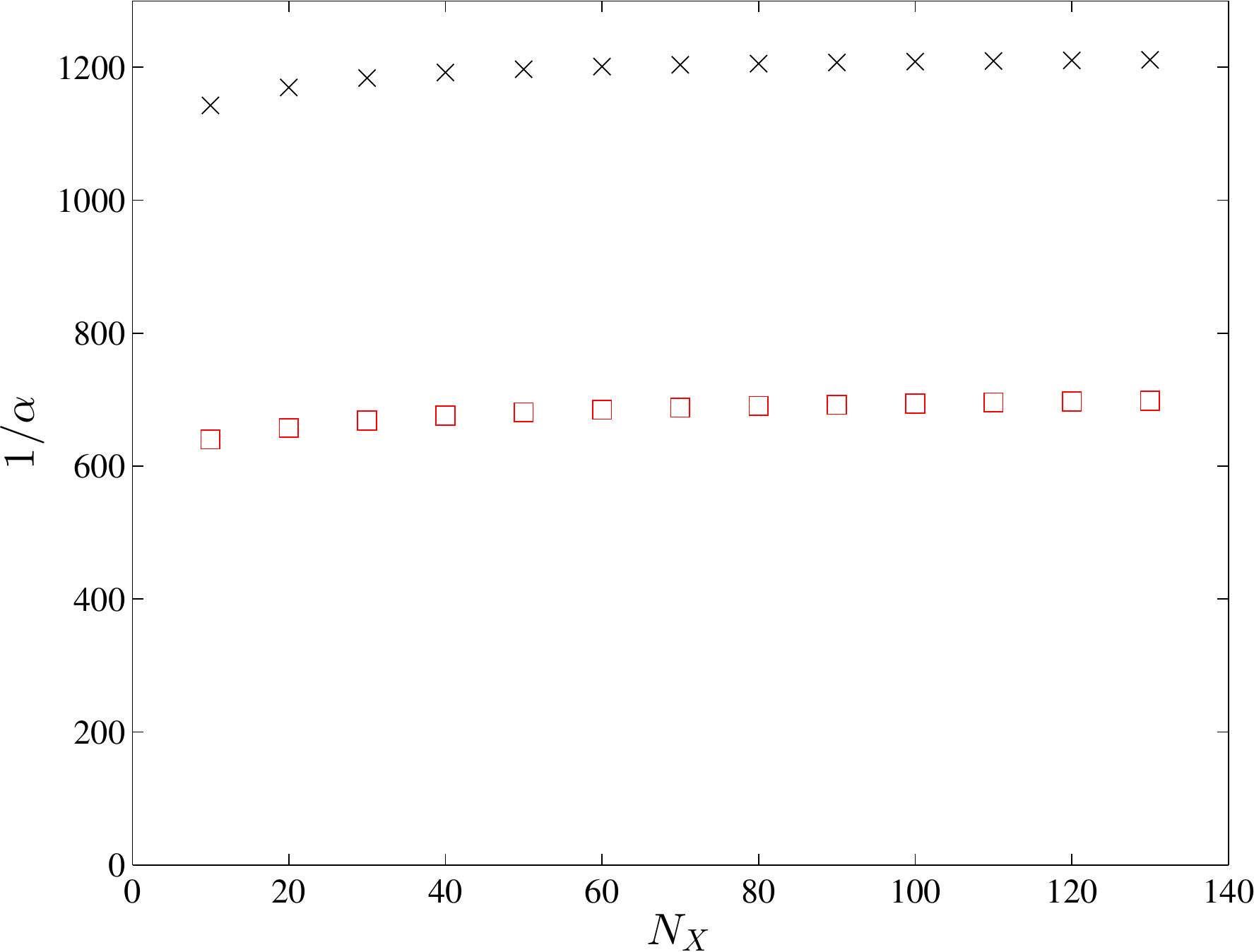}
  \caption{The reciprocal of the step size of example (I). 
  ``{\footnotesize $\square$}'' The Lipschitz constant of $\nabla g$ 
  (used in \textsf{APGM-E}); and 
  ``{\footnotesize $\times$}'' its upper bound
  (used in \textsf{APGM-G}). }
  \label{fig.holonomic_eigenvalue}
\end{figure}

\begin{figure}[tp]
  \centering
  \includegraphics[scale=0.40]{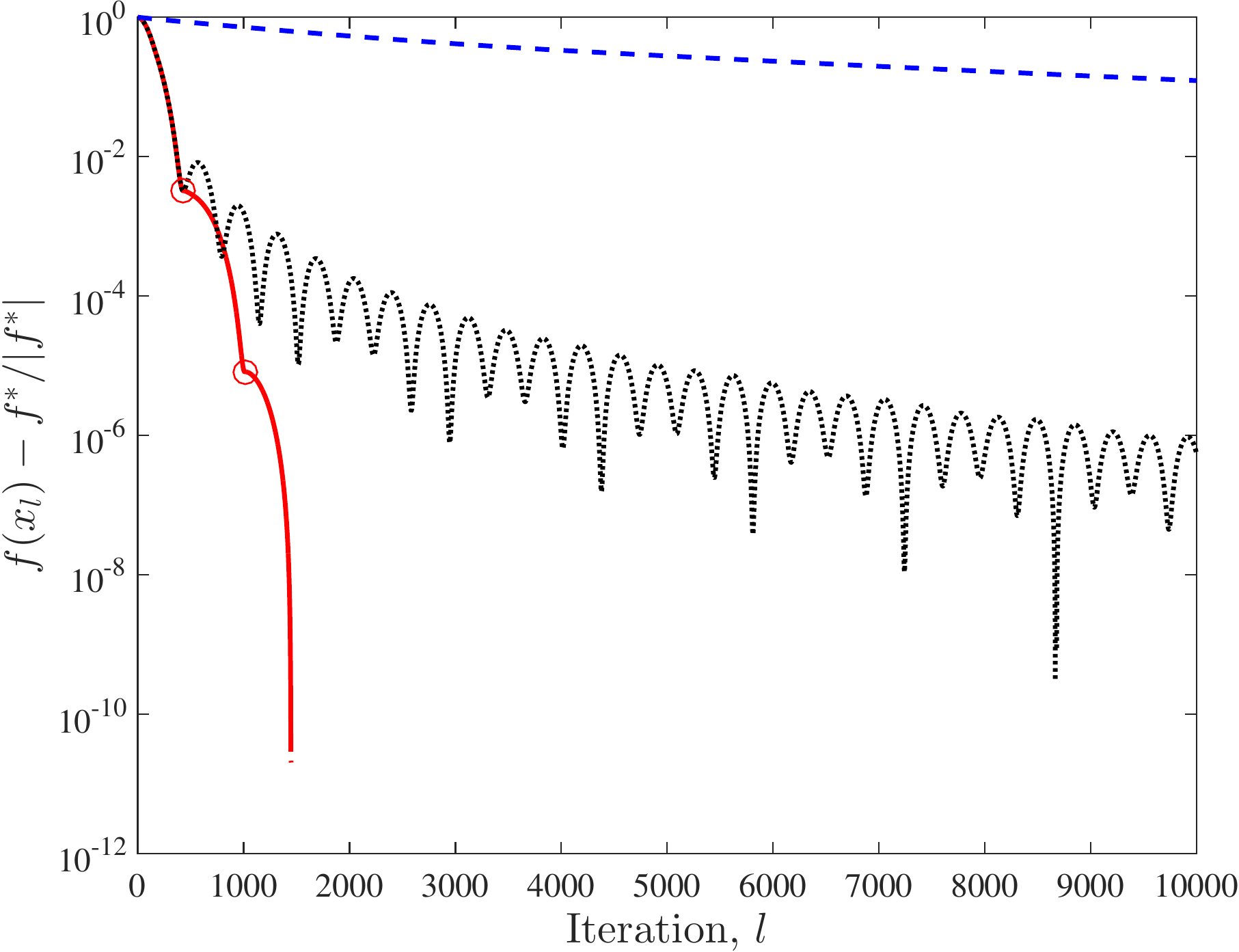}
  \caption{Convergence history of the objective value for 
  $(N_{X},N_{Y})=(40,40)$. 
  ``-----'' \textsf{APGM-E} with restart; 
  ``$\cdot$$\cdot$$\cdot$$\cdot$$\cdot$$\cdot$'' \textsf{APGM-E} without restart; and 
  ``-{\,}-{\,}-'' \textsf{PGM} (\refalg{alg.truss.proximal.1}). }
  \label{fig.restart_comparison}
\end{figure}

In this section we assume a linear hardening model and solve problem 
\eqref{P.truss.3} with \refalg{alg.truss.proximal.3}. 
The initial point chosen at step~\ref{alg.truss.proximal.3.step0} is 
$\bi{v}_{0} = \bi{0}$ and $\bi{p}_{0}=\bi{0}$. 
For comparison, we also solve QP \eqref{P.truss.1} with 
\textsf{QUADPROG}, \textsf{IPOPT}, and \textsf{PATH}. 
It is worth noting that problem \eqref{P.truss.3} has $d+m$ variables, 
while the QP has $d+3m$ variables to be optimized, 
$m$ linear equality constraints, and $2m$ linear inequality constraints. 
As for the external load, $\bi{f}$, a vertical downward force of 
$250/N_{X}$ in $\mathrm{kN}$ is applied at each of the top layer nodes. 

The computational results are listed in \reftab{tab.result.ex.1.1} and 
\reftab{tab.result.ex.1.2}. 
Here, ``iter.''\ means the number of iterations required before 
termination, ``time'' is the total computational time required by each 
algorithm, and ``\texttt{eigs}'' is the computational time required by the 
MATLAB function \texttt{eigs} to compute the maximum eigenvalue of 
$\nabla^{2}g_{1}(\bi{v},\bi{p})$. 
The accuracy of the computed solutions are compared in ``rel.\ diff." 
This reports the value defined by 
$\phi = (\check{f}-f^{*})/f^{*}$, where $f^{*}$ is the objective value 
computed by \textsf{QUADPROG} and $\check{f}$ is the one computed by the 
respective method. 
By definition, $\phi<0$ means that the computed solution has a better 
objective value than the one obtained by a standard QP solver, 
\textsf{QUADPROG}. 
It should be clear that ``time'' of \textsf{APGM-E} includes the 
computational time of \texttt{eigs}. 
For \textsf{PATH}, only the number of major iterations is listed in 
\reftab{tab.result.ex.1.2}, although the number of minor iterations is 
also reported by the solver. 
It is observed from \reftab{tab.result.ex.1.1} and 
\reftab{tab.result.ex.1.2} that the computational time required by 
\textsf{PATH} is extremely large compared with the other four methods. 

\reffig{fig.holonomic_time} shows the computational time of 
\textsf{APGM-E}, \textsf{APGM-G}, \textsf{QUADPROG}, and \textsf{IPOPT}. 
The computational time required by \textsf{IPOPT} is very large compared 
with the other three methods. 
\textsf{QUADPROG} spent two or three times larger time than 
\textsf{APGM-E} and \textsf{APGM-G}. 
\textsf{APGM-E} and \textsf{APGM-G} are comparable from the view point 
of computational time. 
The number of iterations required by these two methods are shown in 
\reffig{fig.holonomic_iter}. 
It is worth noting that the difference of these two methods is only the 
step size, $\alpha$. 
\reffig{fig.holonomic_eigenvalue} shows $1/\alpha$ of these two methods. 
The step size of \textsf{APGM-E} is about $1.75$ times larger than that 
of \textsf{APGM-G}. 
It is observed in \reffig{fig.holonomic_iter} that, for large-scale 
problems, \textsf{APGM-E} requires less iterations; the number of 
iterations required by \textsf{APGM-G} is about $1.3$ times larger. 
Nevertheless, the total computational time is comparable as seen in 
\reffig{fig.holonomic_time}, because in \textsf{APGM-E} the computation 
of the maximum eigenvalue of the Hessian matrix requires relatively 
large computational time. 
In contrast, the computational time required for computing the Gershgorin 
disc bound is negligible (e.g., less than $0.2\,\mathrm{s}$ for the 
instance with $(N_{X},N_{Y})=(130,130)$). 
It is observed in \reftab{tab.result.ex.1.2} that the number of iterations 
of the interior-point method (\textsf{QUADPROG}) is independent of the 
problem size. 
For solving a large-scale instance, the interior-point method has to 
solve a large-scale system of linear equations to find the search 
direction at every iteration, and this computation dominates the 
computational time. 

\reffig{fig.restart_comparison} reports the convergence history of the 
objective value of \textsf{APGM-E} with respect to the iteration count. 
It also shows the result of \textsf{APGM-E} without restart scheme, and 
that of \refalg{alg.truss.proximal.1} (i.e., a proximal gradient method 
without acceleration). 
It is observed that the acceleration and restart schemes drastically 
speed up the convergence.

\subsubsection{Example (II): piecewise-linear hardening model}
\label{sec.ex.2}

\begin{table}[bp]
  \centering
  \caption{The computational results of example (II).}
  \label{tab.holonomic.piecewise.time}
  \begin{tabular}{lrrrrrrrrr}
    \toprule
    & \multicolumn{4}{c}{\textsf{APGM-E}} 
    & \multicolumn{3}{c}{\textsf{APGM-G}} 
    & \multicolumn{2}{c}{\textsf{QUADPROG}} \\
    \cmidrule(lr){2-5} \cmidrule(lr){6-8} \cmidrule(l){9-10}
    $(N_{X},N_{Y})$ 
    & Iter.\ & Time (s) & \texttt{eigs} (s) & Rel.\ diff.\
    & Iter.\ & Time (s) & Rel.\ diff.\ & Iter.\ & Time (s) \\
    \midrule
    $(10,10)$ 
    & 473 & 0.1 & (0.0) & 1.6e-09
    & 394 & 0.1 & 1.2e-06 & 12 & 0.5 \\
    $(20,20)$ 
    & 1061 & 0.7 & (0.2) & 4.3e-07
    & 1526 & 0.8 & 4.4e-07 & 16 & 3.6 \\
    $(30,30)$ 
    & 2159 & 2.5 & (0.6) & 3.1e-08
    & 3102 & 2.7 & 3.4e-08 & 14 & 10.2 \\
    $(40,40)$ 
    & 6181 & 10.1 & (1.3) & 1.0e-07
    & 8861 & 12.5 & 1.1e-07 & 17 & 26.3 \\
    $(50,50)$ 
    & 7763 & 22.3 & (3.1) & 3.1e-07
    & 8118 & 18.1 & 1.4e-05 & 14 & 47.9 \\
    $(60,60)$ 
    & 12391 & 45.7 & (4.7) & 4.2e-08
    & 12981 & 40.1 & 6.2e-06 & 15 & 91.0 \\
    $(70,70)$ 
    & 13263 & 63.8 & (7.6) & 2.9e-06
    & 18971 & 76.4 & 2.9e-06 & 16 & 154.5 \\
    $(80,80)$ 
    & 24767 & 144.6 & (14.6) & 9.4e-09
    & 25646 & 134.1 & 7.6e-06 & 15 & 244.8 \\
    $(90,90)$ 
    & 23616 & 168.7 & (19.0) & 4.2e-05
    & 33734 & 210.9 & 4.2e-05 & 14 & 371.3 \\
    $(100,100)$ 
    & 30782 & 289.3 & (34.8) & 3.7e-05
    & 43946 & 346.1 & 3.7e-05 & 15 & 566.7 \\
    $(110,110)$ 
    & 39675 & 405.0 & (32.1) & 3.3e-06
    & 56613 & 535.9 & 3.3e-06 & 16 & 796.9 \\
    \bottomrule
  \end{tabular}
\end{table}

\begin{figure}[tp]
  \centering
  \includegraphics[scale=0.40]{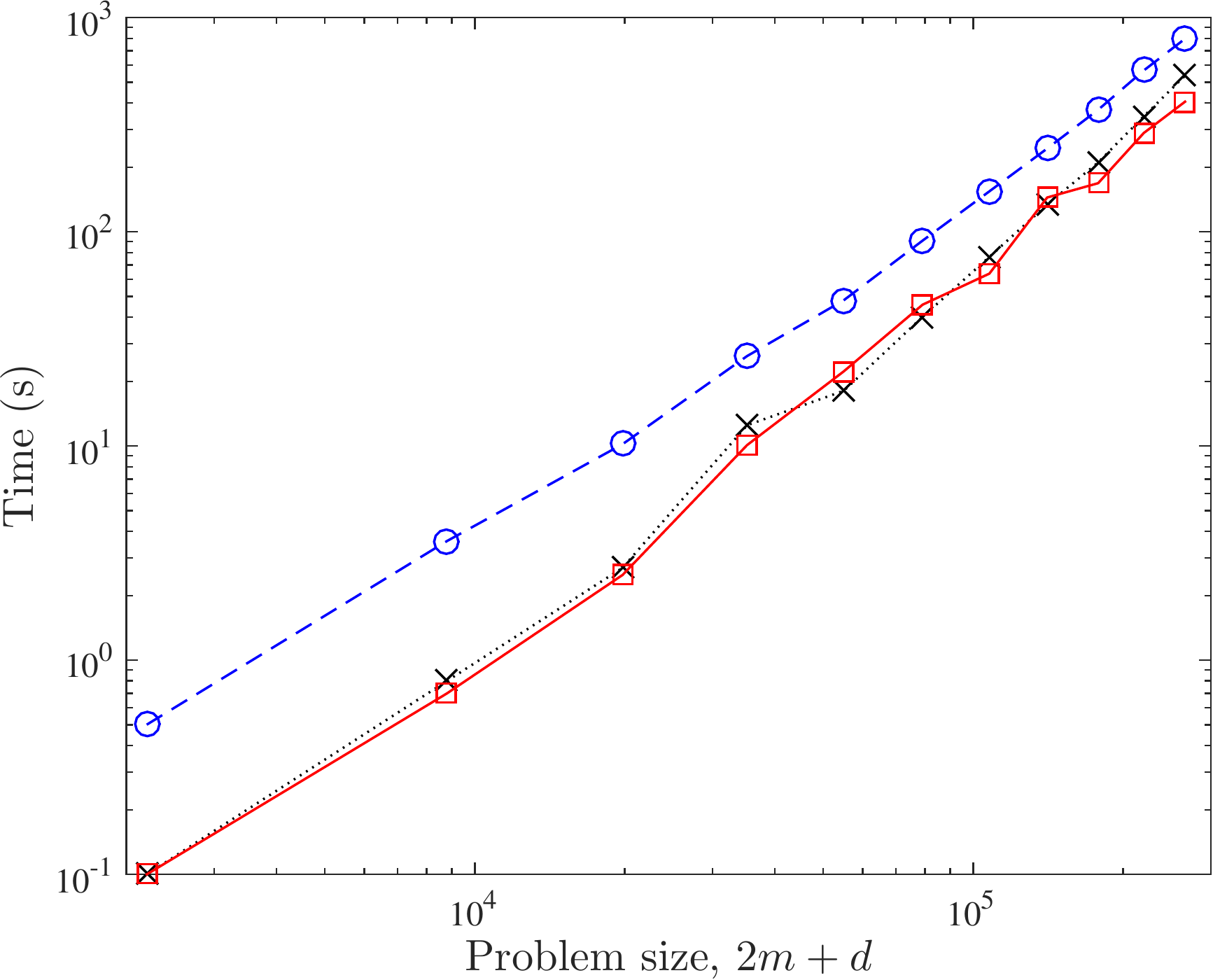}
  \caption{The computational time of example (II). 
  ``{\footnotesize $\square$}'' \textsf{APGM-E}; 
  ``{\footnotesize $\times$}''  \textsf{APGM-G}; and 
  ``{\footnotesize $\circ$}''   \textsf{QUADPROG}. 
  }
  \label{fig.piecewise_time}
\end{figure}

\begin{figure}[tp]
  \centering 
 \includegraphics[scale=0.40]{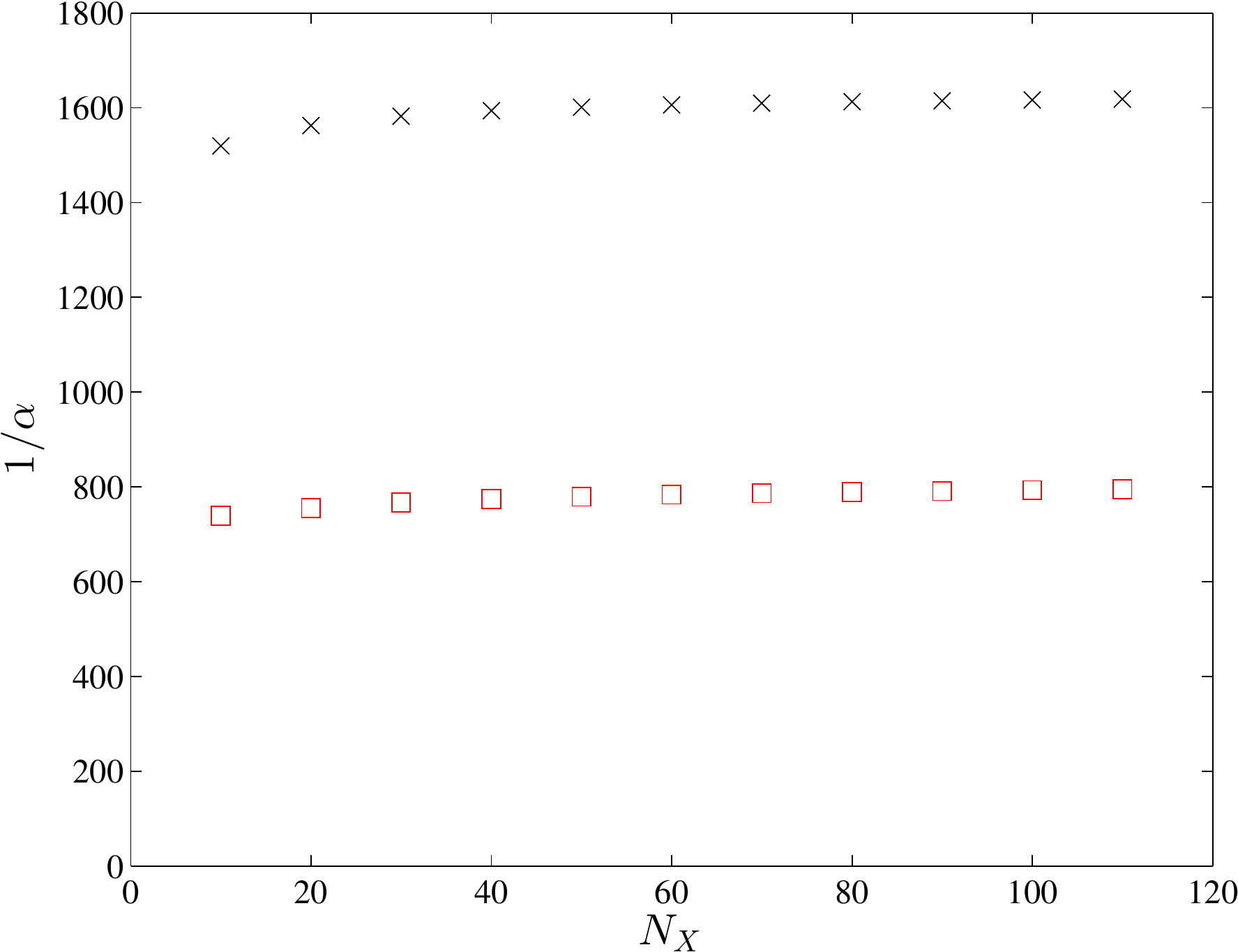}
  \caption{The reciprocal of the step size of example (II). 
  ``{\footnotesize $\square$}'' The Lipschitz constant of $\nabla g$ 
  (used in \textsf{APGM-E}); and 
  ``{\footnotesize $\times$}'' its upper bound
  (used in \textsf{APGM-G}). }
  \label{fig.piecewise_eigenvalue}
\end{figure}

\begin{figure}[tp]
  \centering
  \includegraphics[scale=0.40]{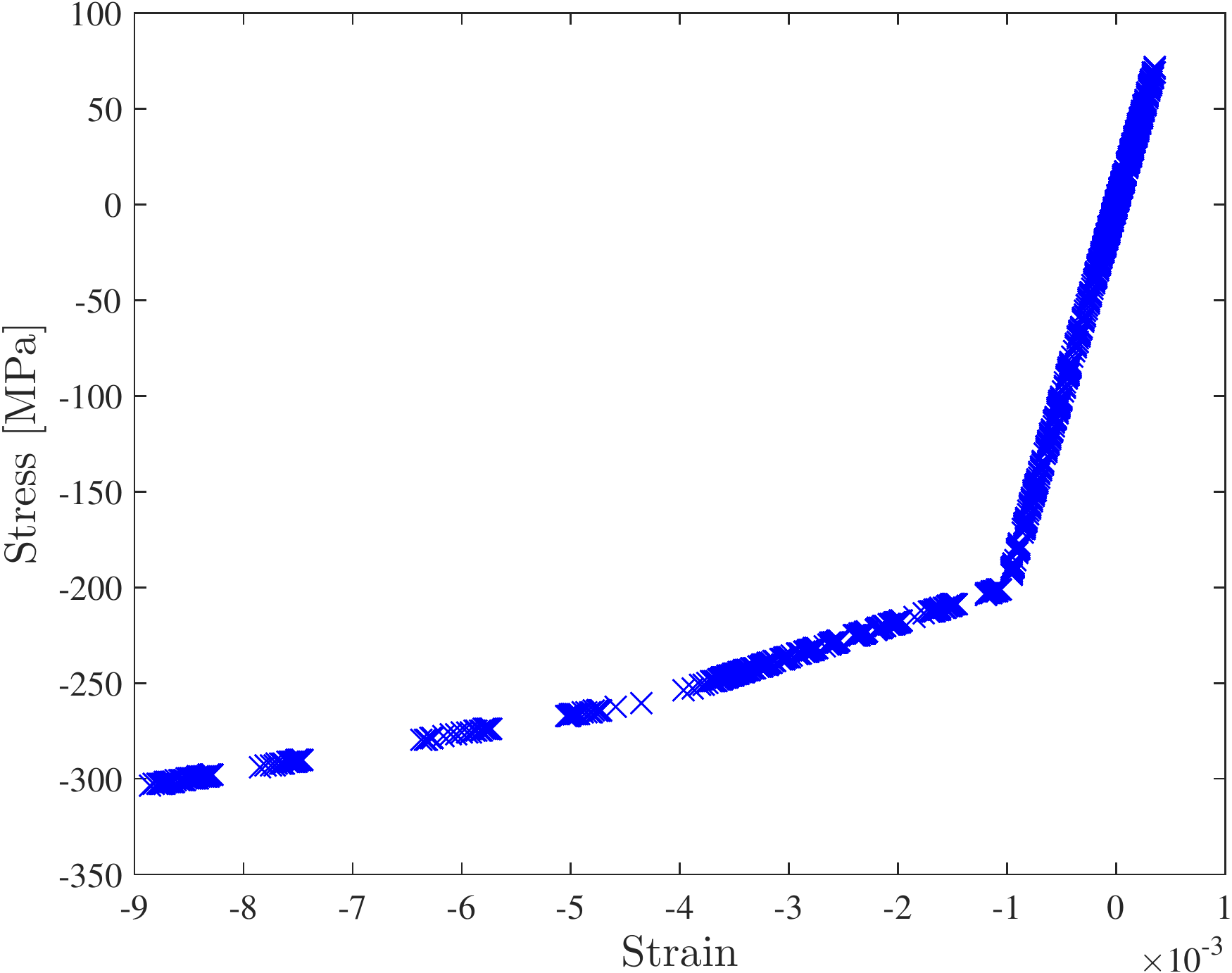}
  \caption{The member stress and strain relations of example (II) with 
  $(N_{X},N_{Y})=(50,50)$. 
  }
  \label{fig.piecewise_stress_strain_50}
\end{figure}

In this section we solve problem instances with a piece-wise linear 
hardening model. 
Specifically, we solve problem \eqref{eq.P.piecewise.1} with 
\refalg{alg.piecewise}. 
The initial point is $\bi{v}_{0}=\bi{0}$ and 
$\bi{p}_{0}=\bi{s}_{0}=\bi{0}$. 
The corresponding QP is problem \eqref{eq.P.piecewise.2}, which is 
solved with an interior-point method. 
It is worth noting that problem \eqref{eq.P.piecewise.1} has $d+2m$ 
variables, while the QP has $d+5m$ variables, 
$m$ linear equality constraints, and $4m$ linear inequality constraints. 
As for the external load, $\bi{f}$, a vertical downward force of 
$200/N_{X}$ (in $\mathrm{kN}$) and a horizontal force of $40/N_{X}$ 
(in $\mathrm{kN}$) in the positive direction of the $Y$-axis are applied 
at each of the top layer nodes. 
The parameters of the hardening model are 
$R^{(0)}_{i}=100\,\mathrm{kN}$, 
$R^{\rr{s}}_{i}=1.3 R^{(0)}_{i}$, 
$h_{i1}=0.1 k_{i}$, and $h_{i2}=0.5 h_{i1}$ $(i=1,\dots,m)$. 

\reftab{tab.holonomic.piecewise.time} lists the computational results. 
The computational time is compared also in \reffig{fig.piecewise_time}. 
\textsf{QUADPROG} spent about twice larger time than \textsf{APGM-E}. 
\textsf{APGM-G} seems to be comparable with \textsf{APGM-E}. 
However, for large instances, \textsf{APGM-G} spent about $1.2$ times 
larger computational time than \textsf{APGM-E}, because 
the number of iterations required by \textsf{APGM-G} is about $1.4$ 
times larger than that of \textsf{APGM-E}. 
\reffig{fig.piecewise_eigenvalue} depicts the maximum eigenvalue of 
$\nabla^{2}g_{1}$ used in \textsf{APGM-E} and its upper bound 
used in \textsf{APGM-G}. 
For all the instances, the upper bound is about twice larger than its 
true value. 
\reffig{fig.piecewise_stress_strain_50} shows the distribution of member 
stress and strain at the equilibrium solution for $(N_{X},N_{Y})=(50,50)$. 
It is observed that the piecewise-linear hardening model is simulated 
correctly.

\subsection{Path-dependent incremental analysis}

In sections~\ref{sec.ex.3} and \ref{sec.ex.4} we perform the 
path-dependent incremental analysis. 
We adopt a combined isotropic and kinematic hardening model studied in 
section~\ref{sec.kinematic}, where $\theta=0.5$. 
Problem \eqref{P.truss.3} is solved with \refalg{alg.truss.proximal.3}. 
We begin by solving the incremental problem for the first loading step 
from an initial point $\bi{v}_{0} = \bi{0}$ and $\bi{p}_{0}=\bi{0}$. 
At each subsequent loading step, we make use of the solution of the 
previous loading step as the initial point. 
The efficiency of this simple warm-start strategy is examined in the 
following examples. 
In the course of incremental analysis, the Hessian matrix of $g_{1}$ is 
independent of the loading step count, $t$. 
Therefore, we need to compute the maximum eigenvalue of the Hessian 
matrix only at the beginning of the analysis at the first loading step. 
Hence, using the maximum eigenvalue outperforms using its 
Gershgorin disc bound.

\subsubsection{Example (III)}
\label{sec.ex.3}

\begin{figure}[tp]
  \centering
  \includegraphics[scale=0.36]{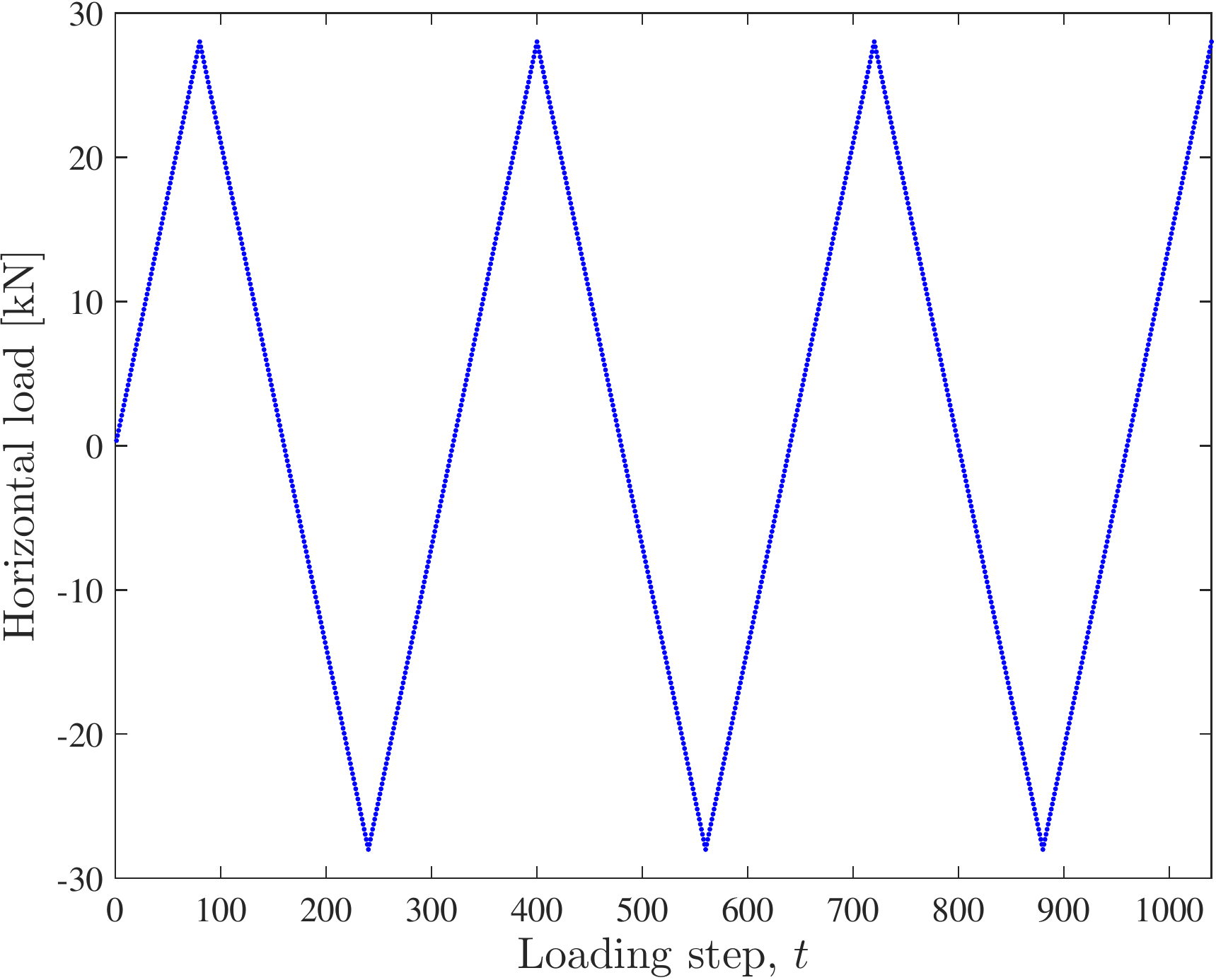}
  \caption{Loading history of example (III). }
  \label{fig.load_x10}
\end{figure}

\begin{figure}[tp]
  \centering
  \subfigure[]{
  \label{fig.displacement_x10}
  \includegraphics[scale=0.32]{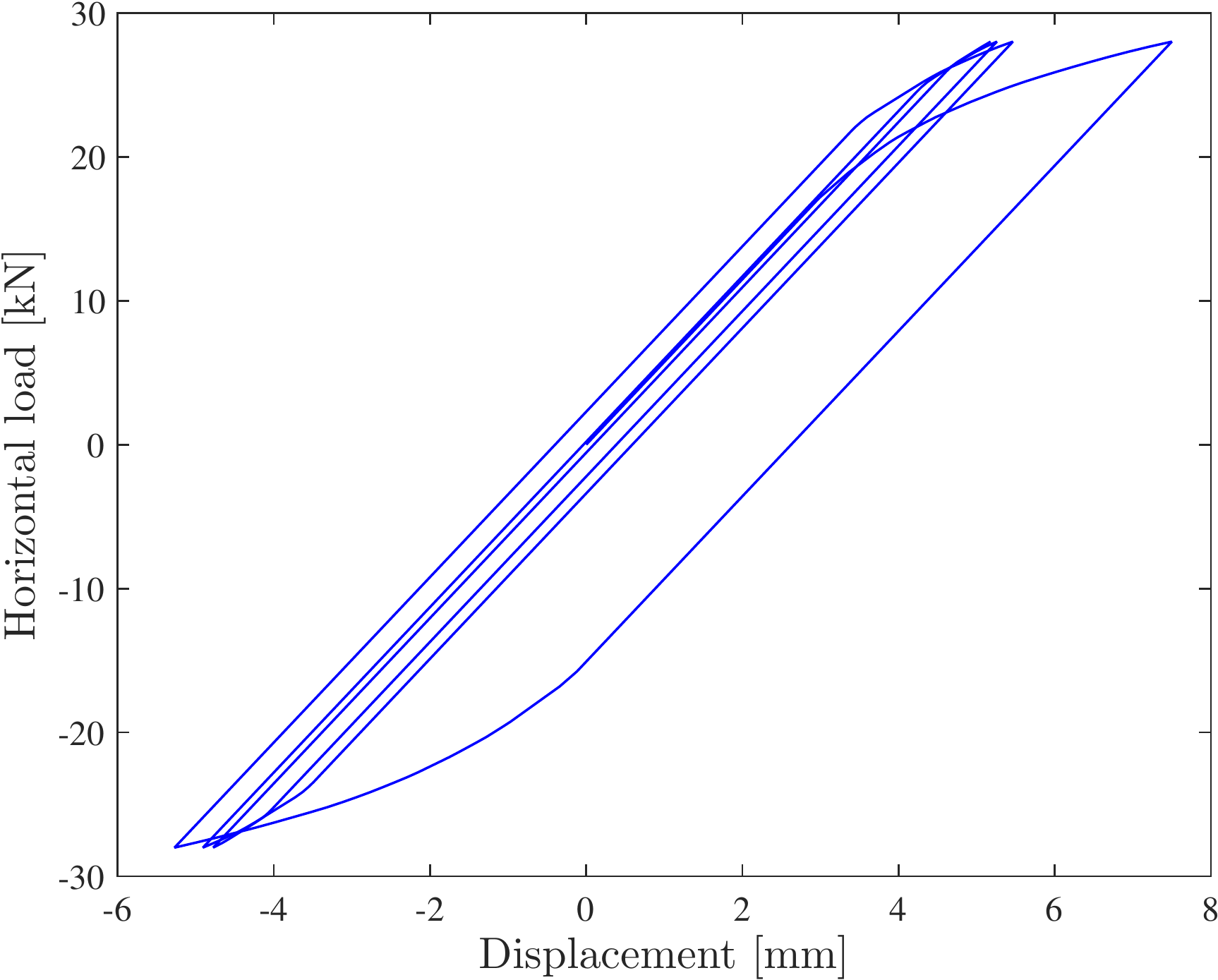}
  }
  \hfill
  \subfigure[]{
  \label{fig.stress_x10}
  \includegraphics[scale=0.32]{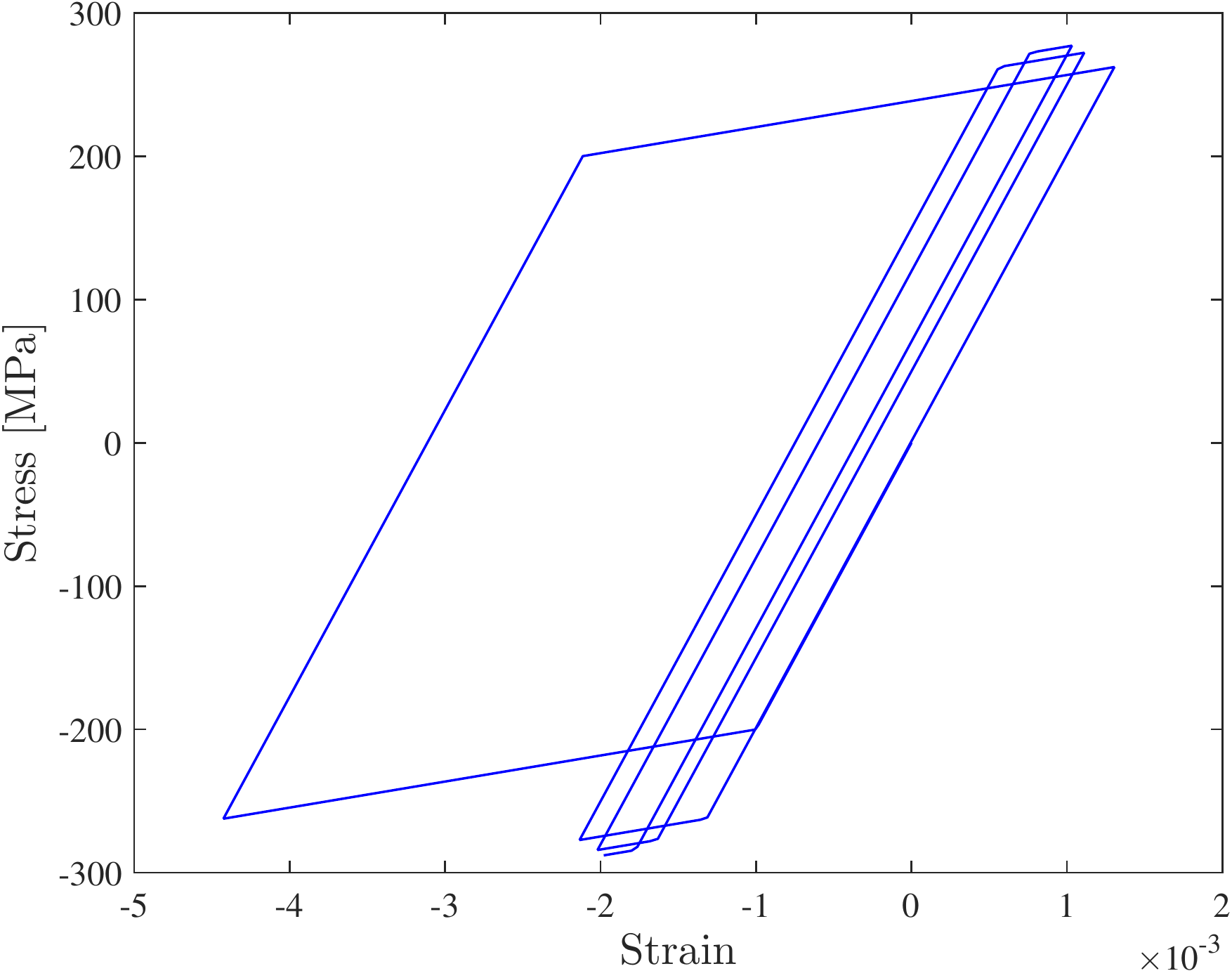}
  }
  \par\bigskip
  \subfigure[]{
  \label{fig.iteration_x10}
  \includegraphics[scale=0.32]{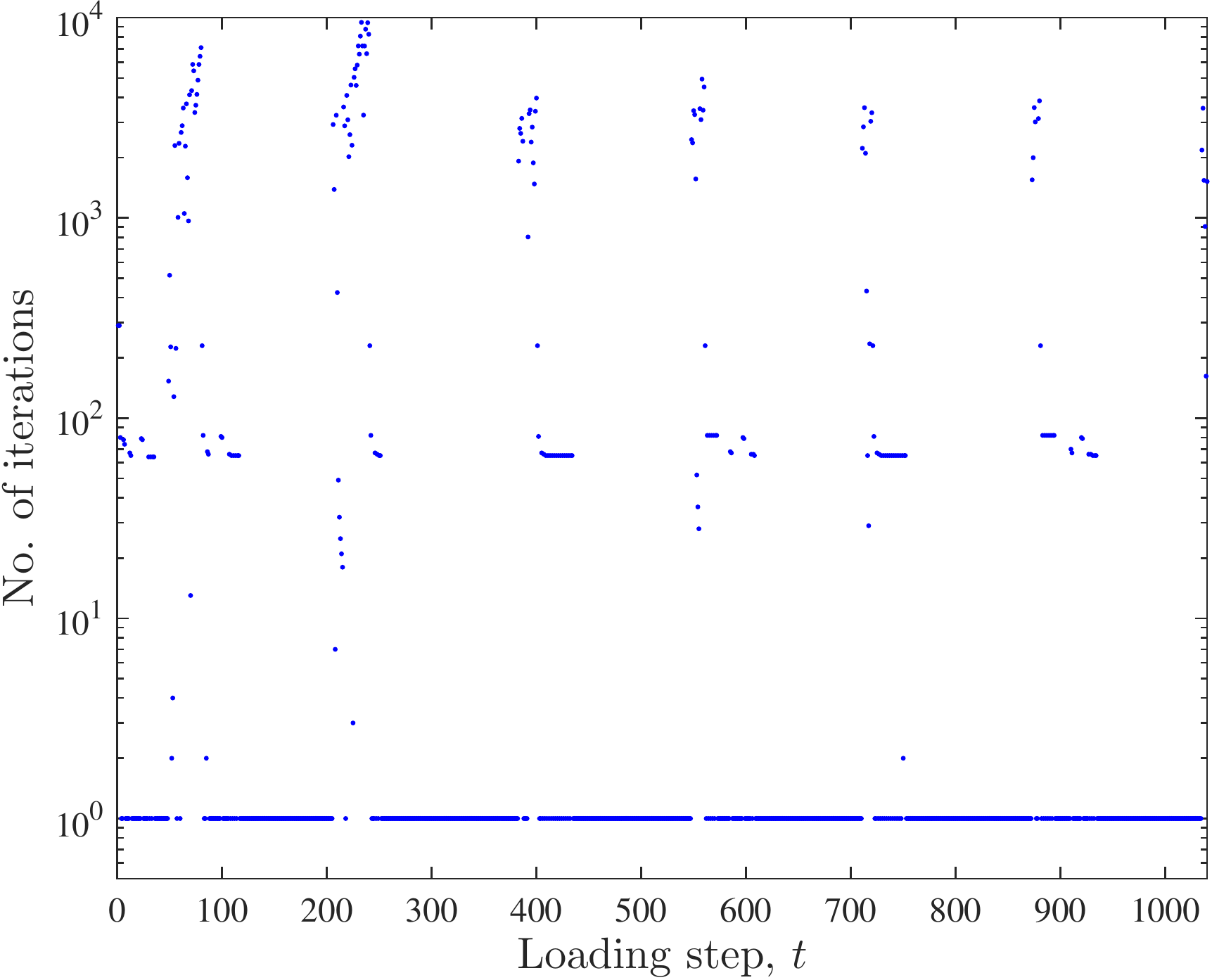}
  }
  \hfill
  \subfigure[]{
  \label{fig.time_x10}
  \includegraphics[scale=0.32]{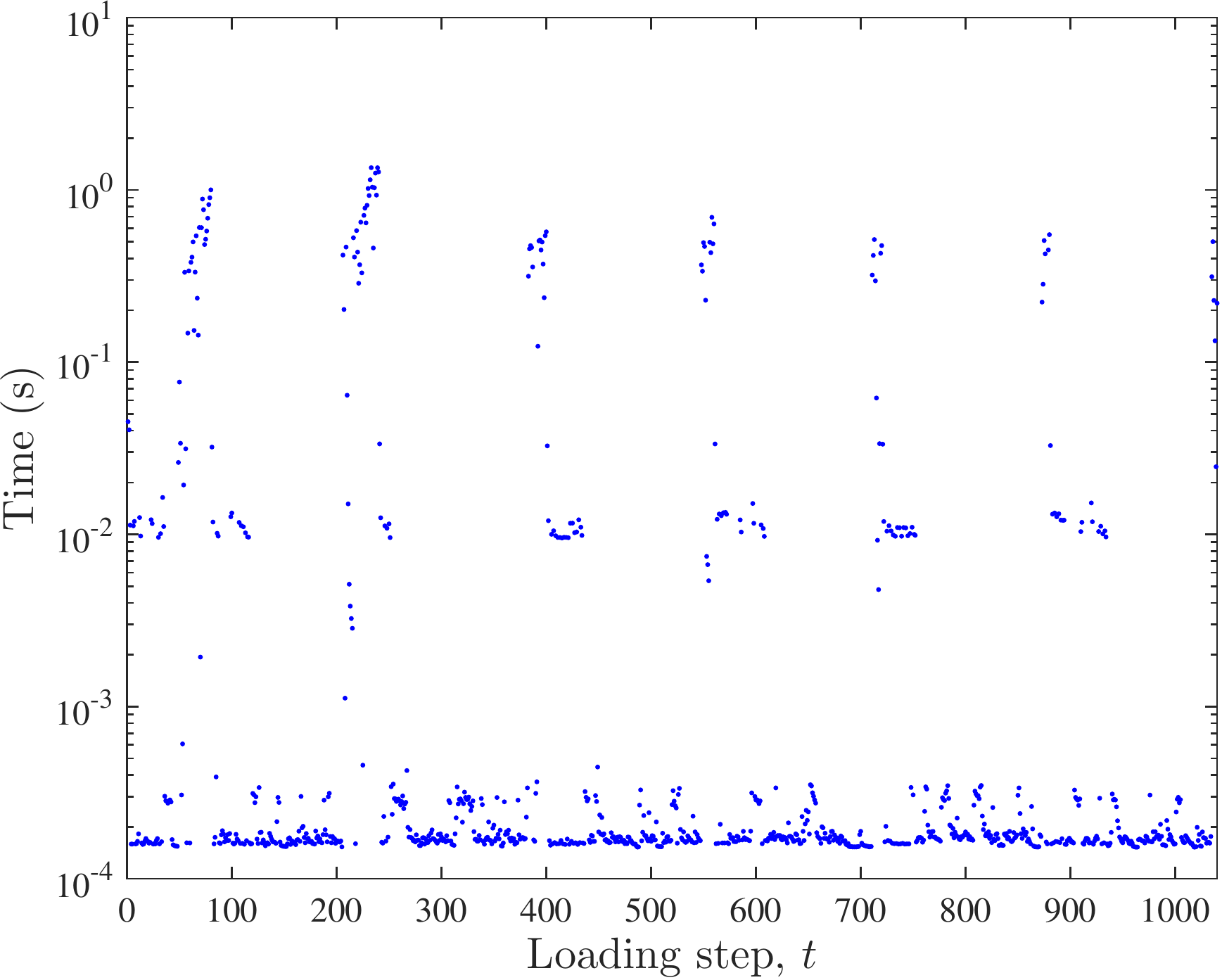}
  }
  \caption{The result of example (III). 
  \subref{fig.displacement_x10} The load versus displacement relation; 
  \subref{fig.stress_x10} the typical stress versus strain relation; 
  \subref{fig.iteration_x10} the number of iterations; and 
  \subref{fig.time_x10} the computational time. 
  }
  \label{fig.incremental_x10}
\end{figure}

In this section we consider $(N_{X},N_{Y})=(10,10)$ for the truss shown 
in \reffig{fig.model_10}. 
A vertical downward force of $5\,\mathrm{kN}$ is applied at each node of 
the top layer. 
Simultaneously, a horizontal force shown in \reffig{fig.load_x10} is 
applied in the positive direction of the $Y$-axis, where $t$ is the 
loading step count. 

\reffig{fig.displacement_x10} shows the load versus the displacement 
relation in the $Y$-direction of the node that is located on 
the $X Z$-plane and on the boundary of the top layer. 
It is observed that the truss gradually shows elastic shakedown, because 
the maximum magnitude of the load is fixed and the isotropic hardening is 
unlimited. 
\reffig{fig.stress_x10} shows the stress--strain relation of a typical 
member. 
The effect of combination of isotropic and kinematic hardening can be 
observed. 
\reffig{fig.iteration_x10} shows the number of iterations required to 
solve the incremental problem at the loading step $t$. 
Similarly, the computational time is shown in \reffig{fig.time_x10}. 
The total computational time was $49.7\,\mathrm{s}$. 
If the incremental solutions both at the $t$th and $(t+1)$th loading 
steps involve no plastic deformation, then these two solutions coincide. 
In such a case, the number of iteration required at the $(t+1)$th 
loading step is negligibly small (and is often one). 
The maximum computational time required for solving one incremental 
problem is $1.87\,\mathrm{s}$.

\subsubsection{Example (IV)}
\label{sec.ex.4}

\begin{figure}[tp]
  \centering
  \includegraphics[scale=0.36]{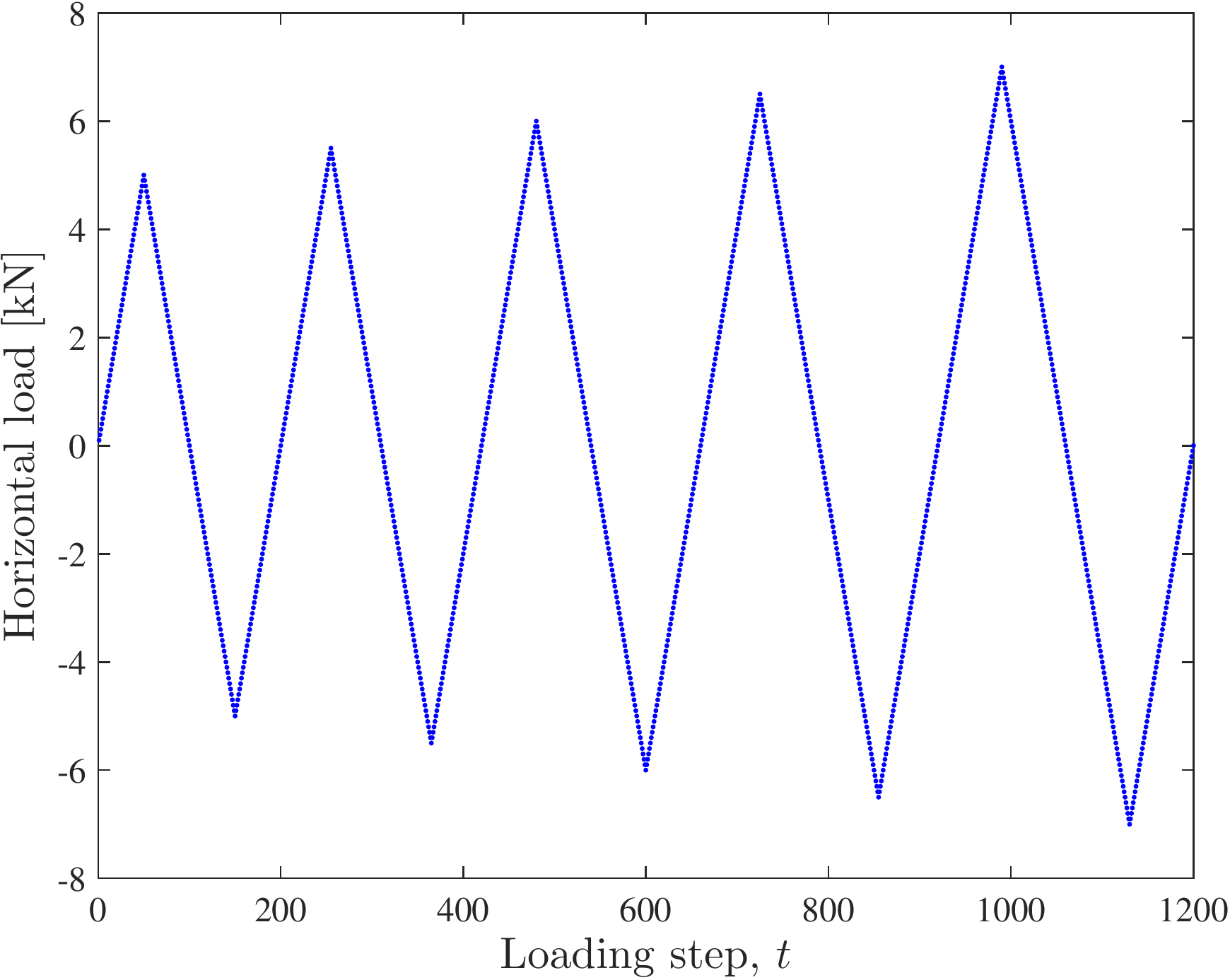}
  \caption{Loading history of example (IV). }
  \label{fig.load_x20}
\end{figure}

\begin{figure}[tp]
  \centering
  \subfigure[]{
  \label{fig.displacement_x20}
  \includegraphics[scale=0.32]{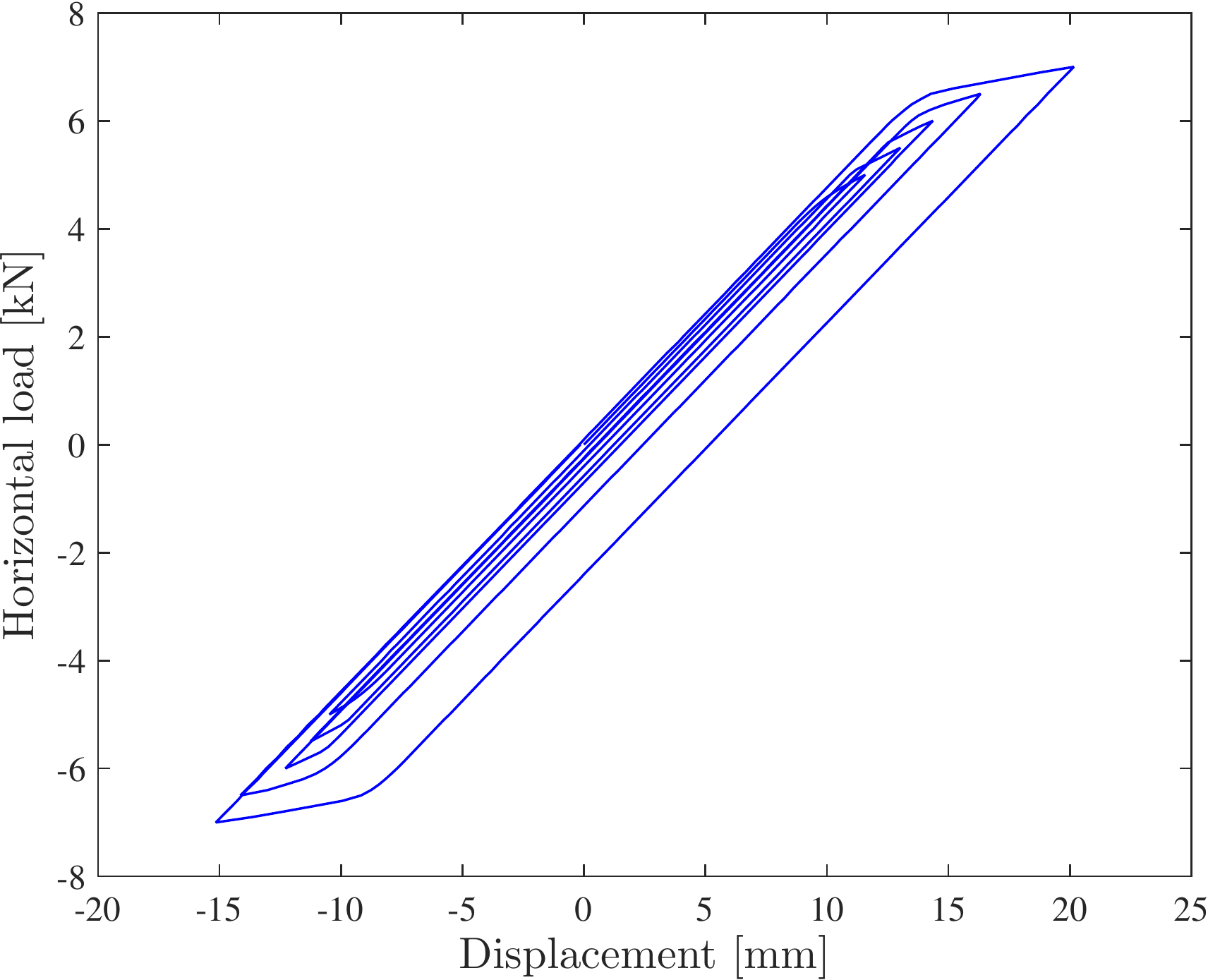}
  }
  \hfill
  \subfigure[]{
  \label{fig.stress_x20}
  \includegraphics[scale=0.32]{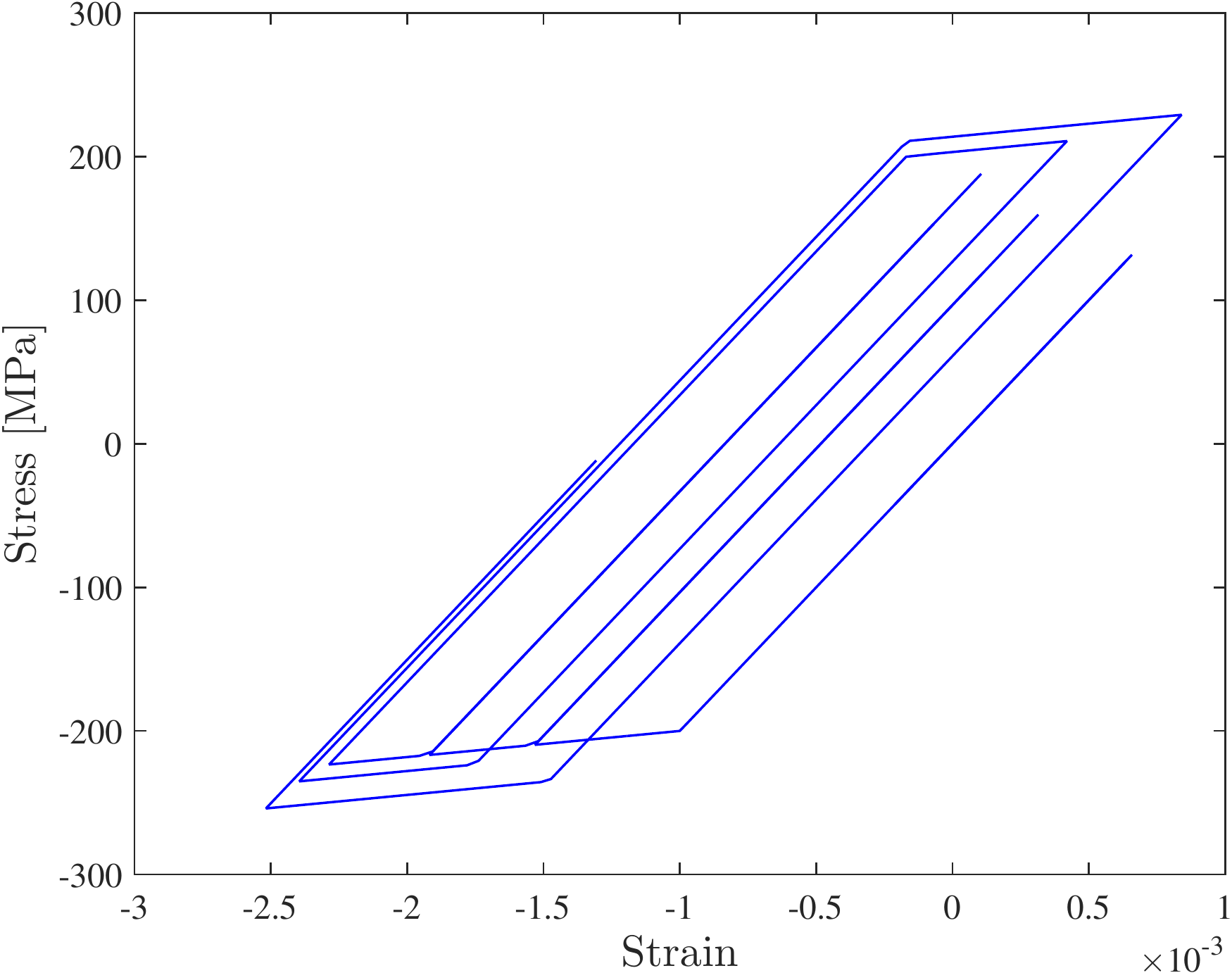}
  }
  \par\bigskip
  \subfigure[]{
  \label{fig.iteration_x20}
  \includegraphics[scale=0.32]{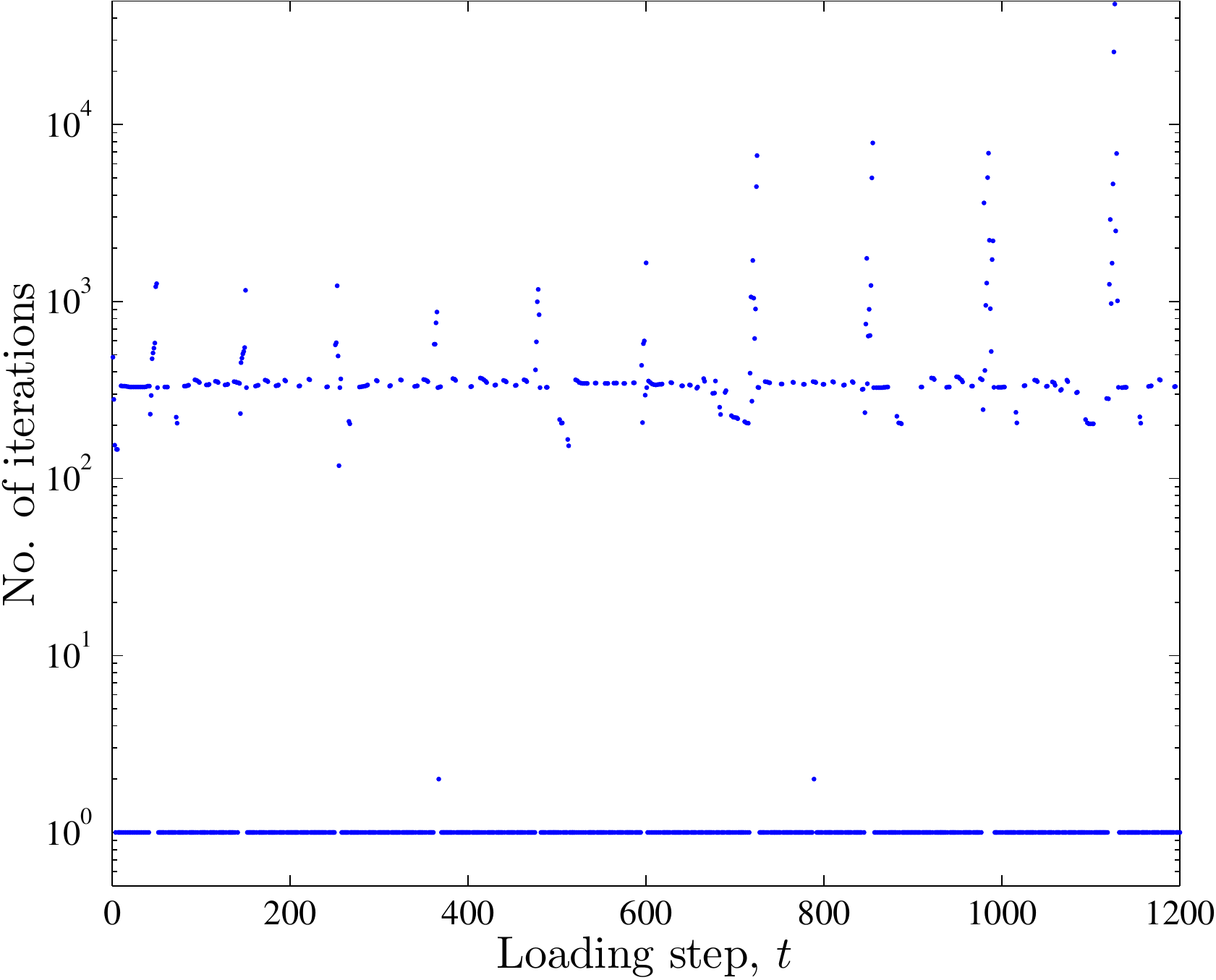}
  }
  \hfill
  \subfigure[]{
  \label{fig.time_x20}
  \includegraphics[scale=0.32]{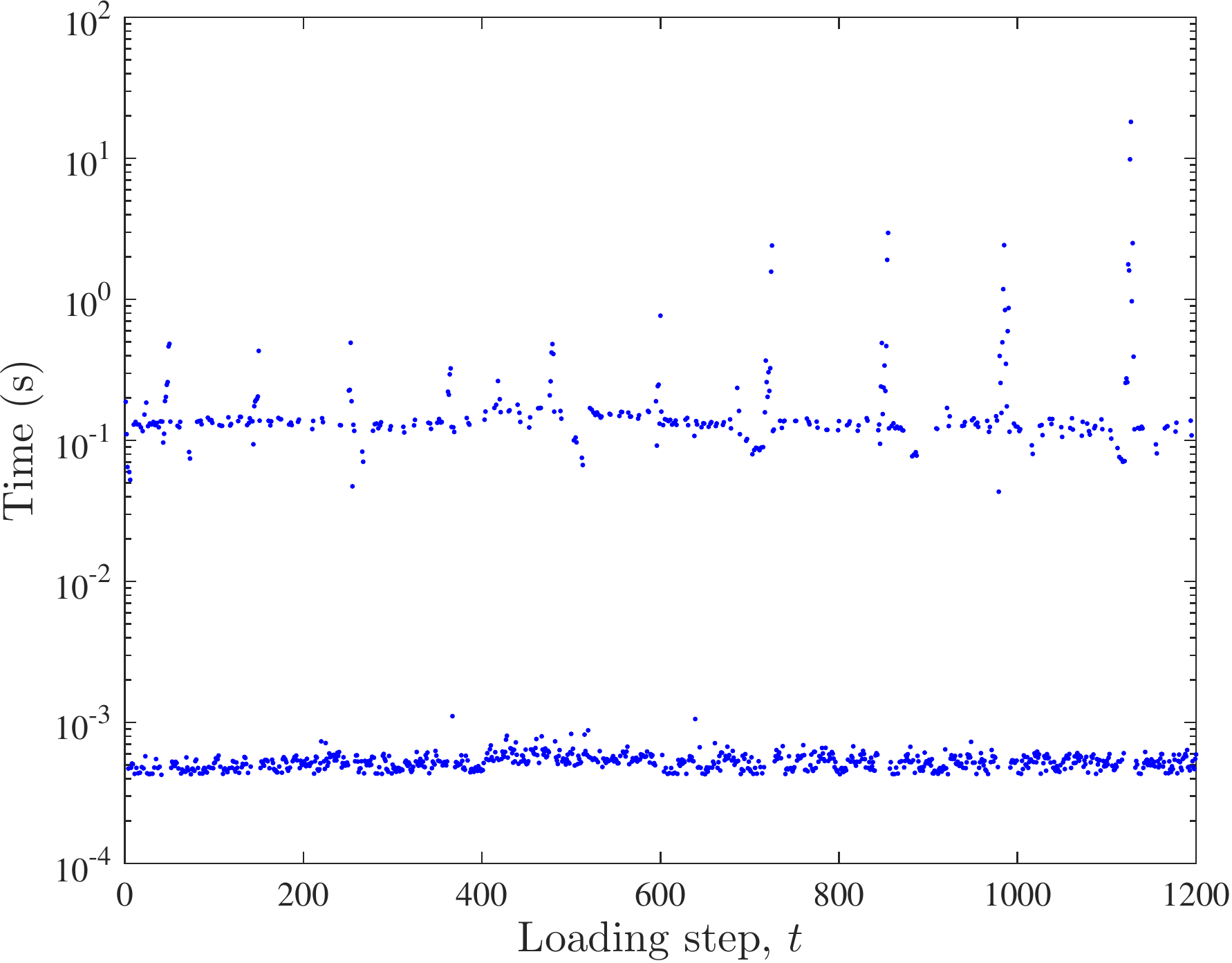}
  }
  \caption{The result of example (IV). 
  \subref{fig.displacement_x20} The load versus displacement relation; 
  \subref{fig.stress_x20} the typical stress versus strain relation; and 
  \subref{fig.iteration_x20} the number of iterations; and  
  \subref{fig.time_x20} the computational time. 
  }
  \label{fig.incremental_x20}
\end{figure}

We next consider a larger instance, $(N_{X},N_{Y})=(20,20)$. 
At each node of the top layer, a vertical downward force of 
$4\,\mathrm{kN}$ is applied. 
Simultaneously, a horizontal force shown in \reffig{fig.load_x20} is 
applied in the $Y$-direction. 

Like \reffig{fig.displacement_x10} in section~\ref{sec.ex.3}, 
\reffig{fig.displacement_x20} shows the load versus the displacement 
relation in the $Y$-direction of the the middle node on the boundary of 
the top layer. 
A typical member stress--strain relation is shown in 
\reffig{fig.stress_x20}. 
The hysteresis loop expands, because the magnitude of the horizontal 
load is gradually increased. 

\reffig{fig.iteration_x20} shows the number of iterations required to 
solve the incremental problem at each loading step. 
There are two cases that required more than $10{,}000$ iterations. 
Loosely speaking, the computational cost increases as the number of 
members with incremental plastic deformations increases. 
The computational time required at each loading step is shown in 
\reffig{fig.time_x20}. 
Since the solution at the previous loading step is used as an initial 
solution for the present loading step, the computational cost becomes 
negligible if no member undergoes plastic deformation. 

The total computational time was $103.2\,\mathrm{s}$. 
Almost all problems were solved within $0.5\,\mathrm{s}$; there exist 
$19$ problems that require more than $0.5\,\mathrm{s}$. 
If we use an interior-point method (\textsf{QUADPROG}) from cold start, 
it is estimated from \reftab{tab.result.ex.1.2} that the total 
computational time might approximately become 
$1.7\,\mathrm{s} \times 1200 = 2040\,\mathrm{s}$. 
The computational time required by the proposed method is much smaller 
than this estimate.

\section{Concluding remarks}
\label{sec.conclude}

In this paper we have presented a fast first-order optimization approach 
to the quasi-static incremental analysis of elastoplastic structures. 
The algorithm is free from numerical solution of linear equations. 
The most expensive computation of the algorithm consists of some 
matrix-vector multiplications with sparse matrices, such as the 
compatibility matrix. 
Also, the presented algorithm does not involve any conditional branching 
procedures stemming from the disjunctive nature of plastic loading and 
elastic unloading processes. 
The algorithm is a version of the accelerated gradient-based methods, 
and converges in potential energy function value as $O(1/k^{2})$, where 
$k$ is the iteration counter. 
More precisely, it is essentially viewed as an application of FISTA, an 
accelerated proximal gradient method for the $\ell_{1}$-regularized 
least-squares problem, to the elastoplastic analysis. 
Owing to these attributes, the algorithm is easy to implement and 
applicable to large-scale problems. 
Indeed, the numerical experiments suggest that the algorithm outperforms 
interior-point methods for convex quadratic programming and nonlinear 
programming. 

In the course of path-dependent quasi-static analysis, we solve 
a series of closely related optimization problems. 
It has been shown that the presented approach can drastically speed up 
by employing a simple warm-start strategy 
that uses the solution at the previous loading step as the initial 
solution for the present loading step. 

It is well known that the incremental problem studied in this paper can 
be recast as a convex quadratic programming problem and a linear 
complementarity problem. 
In contrast, a key to the proposed approach is formulating the 
incremental problem as an unconstrained nonsmooth convex optimization 
problem. 
For simplicity of presentation, in this paper we have restrict ourselves 
to truss structures. 
The presented methodology can be readily applied to other types of 
structures when the yield function can be approximated by a 
piecewise-linear function. 
An example is a frame structure with a piecewise-linear yield condition 
incorporating interaction between the member axial force and end moment. 

This paper has been intended to be the first attempt to apply an 
accelerated gradient-like method to applied mechanics. 
Much remains to be explored. 
For instance, extensions to yield criteria other than the 
piecewise-linear model can be studied; optimization-based approaches to 
such problems can be found in, e.g., 
\citet{KLS07}, \citet{KLSW07}, \citet{TTlS12}, and \citet{YK12}. 
Also, applications of other fast first-order optimization methods can be 
examined. 
Parallelization of the presented method---which is probably quite easy to 
implement because no linear-equations solver is required---has not been 
considered. 
Extension to strain-softening models might be challenging, because it 
requires to deal with a nonconvex objective function, as considered, e.g., 
in \citet{LL15}. 
Recently, it has been discussed that accelerated gradient-like method 
can be viewed as a finite difference approximation of an ordinary 
differential equation \citep{SBC14,KBB15}. 
With reference to these results, the physical interpretation of the 
method presented in this paper might be analyzed. 
Furthermore, besides problems in plasticity theory, extensions to 
complementarity problems arising in diverse fields of nonsmooth 
mechanics can be considered. 
Possible examples include cable networks \citep{KOI02}, 
static and dynamic contact problems \citep{AB08,Wri06}, 
masonry structures \citep{Kan11}, etc. 

More than 30 years ago Giulio Maier wrote \citep{Mai84}: 
  ``{\em Why nonlinear boundary value problems, such as incremental 
  elastoplastic analysis, are routinely solved in several areas of 
  engineering practice fully ignoring the fact that they can be cast in 
  the form of nonlinear or quadratic programs? 
  Obviously, the popular, often merely heuristic, solution schemes 
  resting on iterated use of linear solvers are favoured by the fact 
  that they gradually evolved from the enormous amount of experience 
  accumulated in linear elastic analysis. But their intrinsic 
  superiority over mathematical programming approaches is doubtful, and 
  by no means ensured, in several situations.}'' 
Until today, however, mathematical programming approaches have not been 
used very widely by practitioners. 
The approach presented in this paper has solid background of 
mathematical programming, while computation can be performed 
without knowledge of optimization. 
It might possibly encourage widespread use of various mathematical 
programming approaches to computational and applied mechanics. 


\section*{Acknowledgments}

The author is grateful to Wataru Shimizu for fruitful discussions. 
This work is partially supported by 
JSPS KAKENHI (C) 26420545 
and 
(C) 15KT0109.

%
%
%
%

\appendix

\section{SOCP formulation of problem \eqref{P.truss.1}}
\label{sec.SOCP}

In this section, we explain how problem \eqref{P.truss.1} is recast as a 
second-order cone programming (SOCP) problem. 

The second-order cone in $\Re^{n}$ is defined by 
\begin{align*}
  \LC^{n} = \Bigl\{
  (x_{0},x_{1},\dots,x_{n-1})^{\top} \in \Re^{n} 
  \Bigm|
  x_{0} \ge \sqrt{x_{1}^{2}+\dots+x_{n-1}^{2}} 
  \Bigr\} . 
\end{align*}
SOCP is a minimization (or maximization) of a linear objective function 
under some second-order cone constraints and affine constraints. 

The inequality constraints in \eqref{P.truss.1.3} can be written as 
second-order cone constraints as 
\begin{align*}
  \begin{bmatrix}
    \Delta\gamma_{i} \\  \Delta c_{\rr{p}i} \\
  \end{bmatrix}
  \in \LC^{2} , 
  \quad  i=1,\dots,m . 
\end{align*}
The constraints in \eqref{P.truss.1.2} are affine (i.e., linear equality) 
constraints. 
To convert the objective function to a linear one, we introduce 
auxiliary variables, $\xi \in \Re$ and $\zeta \in \Re$, that serve as 
upper bounds for the quadratic terms in \eqref{P.truss.1.1}. 
Namely, we consider the following constraints: 
\begin{align}
  \xi  &\ge \sum_{i=1}^{m} 
  \frac{1}{2} k_{i} \Delta c_{\rr{e}i}^{2} , 
  \label{eq.reduction.SOCP.1} \\
  \zeta  &\ge \sum_{i=1}^{m} 
  \frac{1}{2} h_{i} \Delta \gamma_{i}^{2} . 
  \label{eq.reduction.SOCP.2}
\end{align}
The convex quadratic inequality constraint in \eqref{eq.reduction.SOCP.1} 
can be rewritten equivalently as \citep{BtN01} 
\begin{align*}
  \xi + 1 \ge 
  \begin{Vmatrix}
    \begin{bmatrix}
      \xi - 1 \\
      \sqrt{2 k_{1}} \Delta c_{\rr{e}1} \\
      \vdots \\
      \sqrt{2 k_{m}} \Delta c_{\rr{e}m} \\
    \end{bmatrix}
  \end{Vmatrix}
  . 
\end{align*}
This is a second-order cone constraint. 
Constraint \eqref{eq.reduction.SOCP.2} can be rewritten in the same 
manner. 

The upshot is that problem \eqref{P.truss.1} can be converted to the 
following SOCP problem:%
\footnote{
Conversion to SOCP is not unique. 
}
  \begin{alignat*}{3}
    & \displaystyle 
    \text{Minimize}
    &{\quad}& \displaystyle
    \sum_{i=1}^{m} q^{(t)}_{i} \Delta c_{\rr{e}_{i}} 
    + \xi 
    + \sum_{i=1}^{m} R^{(t)}_{i} \Delta\gamma_{i} 
    + \zeta 
    - \bi{f}^{\top} \Delta\bi{u} \\
    & \ST && \displaystyle
    \Delta c_{\rr{e}i} + \Delta c_{\rr{p}i} 
    = \bi{b}_{i}^{\top} \Delta\bi{u} , 
    \quad  i=1,\dots,m, \\
    & && \displaystyle
    \begin{bmatrix}
      \Delta\gamma_{i} \\  \Delta c_{\rr{p}i} \\
    \end{bmatrix}
    \in \LC^{2} , 
    \quad  i=1,\dots,m , \\
    & && \displaystyle
    \begin{bmatrix}
      \xi + 1 \\
      \xi - 1 \\
      \sqrt{2 k_{1}} \Delta c_{\rr{e}1} \\
      \vdots \\
      \sqrt{2 k_{m}} \Delta c_{\rr{e}m} \\
    \end{bmatrix}
    \in \LC^{m+1} , 
    \quad
    \begin{bmatrix}
      \zeta + 1 \\
      \zeta - 1 \\
      \sqrt{2 h_{1}} \Delta \gamma_{1} \\
      \vdots \\
      \sqrt{2 h_{m}} \Delta \gamma_{m} \\
    \end{bmatrix}
    \in \LC^{m+1} . 
  \end{alignat*}%
Here, $\Delta c_{\rr{e}_{1}},\dots,\Delta c_{\rr{e}_{m}}$, 
$\xi$, $\Delta\gamma_{1},\dots,\Delta\gamma_{m}$, $\zeta$, and 
$\Delta\bi{u}$ are variables to be optimized.

\section{Equivalence of \eqref{eq.optimality.proximal.3} and \eqref{eq.optimality.proximal.6}}
\label{sec.fixed.point}
As one of fundamental properties of the proximal mapping, 
we can show, for any $\alpha > 0$, that $\bi{p} \in \Re^{m}$ satisfies 
\begin{align}
  \bi{0} &\in 
  \nabla_{\bi{p}}g_{1}(\bi{v},\bi{p}) 
  + \partial g_{2}(\bi{p})  
  \label{eq.optimality.proof.1}
\end{align}
if and only if it satisfies 
\begin{align}
  \bi{p}
  &= \bs{prox}_{\alpha g_{2}}
  (\bi{p} - \alpha\nabla_{\bi{p}} g_{1}(\bi{v},\bi{p})) . 
  \label{eq.optimality.proof.2}
\end{align}
See, e.g., \citet{PB14}. 
For the reader's convenience, essentials of the proof are repeated here. 

Suppose that $\bi{p}$ satisfies \eqref{eq.optimality.proof.1}. 
This is equivalent to 
\begin{align}
  \bi{0} &\in 
  \alpha \nabla_{\bi{p}}g_{1}(\bi{v},\bi{p}) 
  + \alpha \partial g_{2}(\bi{p})  \notag\\
  &= 
  \alpha \nabla_{\bi{p}}g_{1}(\bi{v},\bi{p}) 
  - \bi{p} + \bi{p} 
  + \alpha \partial g_{2}(\bi{p}) . 
  \label{eq.optimality.proof.3}
\end{align}
Let $\bi{s} := \bi{p} - \alpha \nabla_{\bi{p}}g_{1}(\bi{v},\bi{p})$ for 
notational simplicity. 
Then \eqref{eq.optimality.proof.3} is rewritten as 
\begin{align*}
  \bi{0} &\in 
  \alpha \partial g_{2}(\bi{p}) + (\bi{p} - \bi{s}) , 
\end{align*}
which is equivalent to 
\begin{align}
  \bi{p} = \argmin_{\bi{z}} 
  \Bigl\{ \alpha g_{2}(\bi{z}) 
  + \frac{1}{2} \| \bi{z} - \bi{s} \|^{2} \Bigr\} . 
  \label{eq.optimality.proof.4}
\end{align}
By definition, \eqref{eq.optimality.proof.4} is equivalent to 
\eqref{eq.optimality.proof.2}.

\section{Algorithm for piecewise-linear hardening}
\label{sec.algorithm_piece}
We begin with computation of the gradient of $g_{1}$ defined by 
\eqref{eq.P.piecewise.g.1}. 
In a manner similar to section~\ref{sec.accelerated}, 
it is convenient to define $\bi{e} \in \Re^{m}$ by 
\begin{align*}
  \bi{e} 
  = \ti{B} \bi{v} - \bi{p} - \bi{s} , 
\end{align*}
which corresponds to the vector of incremental elastic elongation, 
$\Delta\bi{c}_{\rr{e}}$, in problem \eqref{eq.P.piecewise.2}. 
Then the gradient of $g_{1}$ can be calculated as 
\begin{align*}
  \nabla_{\bi{v}}g_{1}(\bi{v},\bi{p},\bi{s}) 
  &= \ti{B}^{\top} (\diag(\bi{k}) \bi{e} + \bi{q}^{(t)}) 
  - \bi{f} ,  
  \\
  \nabla_{\bi{p}}g_{1}(\bi{v},\bi{p},\bi{s}) 
  &= \diag(\bi{h}_{1}) \bi{p}  - \diag(\bi{k}) \bi{e} - \bi{q}^{(t)} , 
  \\
  \nabla_{\bi{s}}g_{1}(\bi{v},\bi{p},\bi{s}) 
  &= \diag(\bi{\eta}) \bi{s}  - \diag(\bi{k}) \bi{e} - \bi{q}^{(t)} , 
\end{align*}
where 
\begin{align*}
  \nabla_{\bi{v}}g_{1} 
  = \pdif{g_{1}}{\bi{v}} , 
  \quad
  \nabla_{\bi{p}}g_{1}
  = \pdif{g_{1}}{\bi{p}} , 
  \quad
  \nabla_{\bi{s}}g_{1}
  = \pdif{g_{1}}{\bi{s}} . 
\end{align*}
Moreover, the Hessian matrix of $g_{1}$ is 
written as 
\begin{align}
  \nabla^{2} g_{1}(\bi{v},\bi{p},\bi{s}) 
  &= 
  \left[\begin{array}{@{}c|c|c@{\,}}
   \displaystyle
     \ti{B}^{\top} & \ti{O} & \ti{O} \\
     \hline
     -\ti{I} & \ti{I} & \ti{O} \\
     \hline
     -\ti{I} & \ti{O} & \ti{I} \\
  \end{array}
  \right] 
  \left[\begin{array}{@{}c|c|c@{\,}}
   \displaystyle
     \diag(\bi{k}) & \ti{O} & \ti{O} \\
     \hline
     \ti{O} & \diag(\bi{h}_{1}) & \ti{O} \\
     \hline
     \ti{O} & \ti{O} & \diag(\bi{\eta}) \\
  \end{array}
  \right] 
  \left[\begin{array}{@{}c|c|c@{\,}}
   \displaystyle
     \ti{B} & -\ti{I} & -\ti{I} \\
     \hline
     \ti{O} & \ti{I} & \ti{O} \\
     \hline
     \ti{O} & \ti{O} & \ti{I} \\
  \end{array}
  \right] . 
  \label{eq.nabla.2.g1.piecewise}
\end{align}
Since $k_{i}>0$, $h_{i1}>0$, $\eta_{i}>0$ $(i=1,\dots,m)$ and $\ti{B}$ 
is of row full rank for a kinematically determinate truss, 
$\nabla^{2} g_{1}(\bi{v},\bi{p},\bi{s})$ is positive definite. 
In a manner similar to section~\ref{sec.accelerated}, the proximal 
mapping of $\alpha g_{2}$ with $\alpha > 0$ can be computed as 
\begin{align*}
  \bs{prox}_{\alpha g_{2}}(\bi{w},\bi{z}) = 
  \begin{bmatrix}
    \diag(\sign(\bi{w})) 
    \max \{ |\bi{w}| - \alpha \bi{R}^{(t)}, \bi{0} \} \\
    \diag(\sign(\bi{z})) 
    \max \{ |\bi{z}| - \alpha \bi{R}^{\rr{s}}, \bi{0} \} \\
  \end{bmatrix}
  . 
\end{align*}

We are now in position to describe an accelerated proximal gradient 
method for solving problem \eqref{eq.P.piecewise.1}. 
\begin{algorithm}\label{alg.piecewise}
  \hspace{1em}
  \begin{algstep}
    \setcounter{alnum}{-1}
    \item 
    Let $L$ denote the maximum eigenvalue of 
    $\nabla^{2} g_{1}(\bi{v},\bi{p},\bi{s})$ in \eqref{eq.nabla.2.g1.piecewise}. 
    Choose $\bi{v}_{0} \in \Re^{d}$, $\bi{p}_{0} \in \Re^{m}$, 
    $\bi{s}_{0} \in \Re^{m}$, $\alpha \in ]0,1/L]$, and termination 
    tolerance $\epsilon >0$. 
    Set $l:=1$, $\bi{\mu}_{1} := \bi{v}_{0}$, 
    $\bi{\rho}_{1} := \bi{p}_{0}$, 
    $\bi{\sigma}_{1} := \bi{s}_{0}$, and 
    $\tau_{1}:=1$. 
    
    \item \label{alg.piecewise.step.1}
    Let 
    \begin{align*}
      \bi{\varepsilon}_{l} 
      &:= \ti{B} \bi{\mu}_{l} - \bi{\rho}_{l} - \bi{\sigma}_{l} , \\
      \bi{v}_{l} 
      &:= \bi{\mu}_{l} - \alpha 
      [ \ti{B}^{\top} (\diag(\bi{k}) \bi{\varepsilon}_{l} + \bi{q}^{(t)}) 
      - \bi{f} ] , \\
      \bi{w}_{l}
      &:= \bi{\rho}_{l} - \alpha 
      (\diag(\bi{h}_{1}) \bi{\rho}_{l} - \diag(\bi{k}) \bi{\varepsilon}_{l} 
      - \bi{q}^{(t)}) , \\
      \bi{p}_{l} 
      &:= \diag(\sign(\bi{w}_{l})) \max \{ |\bi{w}_{l}| - \alpha \bi{R}^{(t)}, \bi{0} \} , \\
      \bi{z}_{l}
      &:= \bi{\sigma}_{l} - \alpha 
      (\diag(\bi{\eta}) \bi{\sigma}_{l} - \diag(\bi{k}) \bi{\varepsilon}_{l} 
      - \bi{q}^{(t)}) , \\
      \bi{s}_{l} 
      &:= \diag(\sign(\bi{z}_{l})) \max \{ |\bi{z}_{l}| - \alpha \bi{R}^{\rr{s}}, \bi{0} \}.
    \end{align*}
    
    \item 
    Let 
    \begin{align*}
      \tau_{l+1}  
      = \frac{1}{2} \Bigl( 1 + \sqrt{1 + 4 \tau_{l}^{2}} \Bigr) . 
    \end{align*}
    
    \item 
    If $g_{1}(\bi{v}_{l},\bi{p}_{l},\bi{s}_{l}) 
    + g_{2}(\bi{p}_{l},\bi{s}_{l}) 
    < g_{1}(\bi{v}_{l-1},\bi{p}_{l-1},\bi{s}_{l-1}) 
    + g_{2}(\bi{p}_{l-1},\bi{s}_{l-1})$, 
    then let 
    \begin{align*}
      \bi{\mu}_{l+1}  
      &:= \bi{v}_{l} 
      + \frac{\tau_{l}-1}{\tau_{l+1}} (\bi{v}_{l} - \bi{v}_{l-1}) , \\
      \bi{\rho}_{l+1}  
      &:= \bi{p}_{l} + \frac{\tau_{l}-1}{\tau_{l+1}} 
      (\bi{p}_{l} - \bi{p}_{l-1}) , \\
      \bi{\sigma}_{l+1}  
      &:= \bi{s}_{l} + \frac{\tau_{l}-1}{\tau_{l+1}} 
      (\bi{s}_{l} - \bi{s}_{l-1}) . 
    \end{align*}
    Otherwise, let 
    $\tau_{l+1}:=1$, $\bi{\mu}_{l+1}:=\bi{v}_{l}$, 
    $\bi{\rho}_{l+1}:=\bi{p}_{l}$, and 
    $\bi{\sigma}_{l+1}:=\bi{s}_{l}$. 
    
    \item
    If $\| (\bi{v}_{l}, \bi{p}_{l}, \bi{s}_{l}) 
    - (\bi{v}_{l-1}, \bi{p}_{l-1}, \bi{s}_{l-1}) \| \le \epsilon$, 
    then terminate. 
    Otherwise, let $l \gets l+1$, and go to 
    step~\ref{alg.piecewise.step.1}. 
  \end{algstep}
\end{algorithm}

At step~\ref{alg.piecewise.step.1} of \refalg{alg.piecewise}, auxiliary 
variables $\bi{\varepsilon}_{l}$, $\bi{w}_{l}$, and $\bi{s}_{l}$ are 
used for convenience of computation. 

%
%

\end{document}